\documentclass[11 pt]{amsart}
\usepackage{mathrsfs}
 \newtheorem{thm}{Theorem}[section]

 \theoremstyle{definition}
 
 \theoremstyle{remark}

 \numberwithin{equation}{section}

\begin{document}
\title[Maps preserving   numerical range of Lie product]{Non-linear maps on self-adjoint operators preserving  numerical radius and numerical range  of Lie product}

\author{Jinchuan Hou}
\address[Jinchuan Hou]{College of
 Mathematics, Taiyuan University of Technology, Taiyuan,
 030024, P. R. China} \email[J. Hou]{jinchuanhou@aliyun.com}

\author{Kan He}
\address[Kan He]{College of
Mathematics, Taiyuan University of Technology, Taiyuan,
 030024, P. R. China} \email[K.
He]{hk19830310@163.com}

\thanks{This work supported by National
Science Foundation of China ( 11171249, 11201329,11271217) and
Program for the Outstanding Innovative Teams of Higher Learning
Institutions of Shanxi.}
 \maketitle
\begin{abstract}

Let $H$ be a  complex separable Hilbert space of dimension $\geq 2$,
${\mathcal B}_s(H)$ the space of all self-adjoint operators on $H$.
We give a complete classification of non-linear surjective maps on
$\mathcal B_s(H)$ preserving respectively numerical radius and
numerical range of Lie product.

\end{abstract}

\section{Introduction}

Let $A$ be a bounded linear operator acting on a complex Hilbert
space $H$. Recall that the
 numerical range of $A$ is the set $W(A)=\{\langle Ax,x\rangle \,|\, x\in H, \|x\|=1 \},$
 and the numerical radius of $A$ is $w(A)=\sup\{|\lambda| \,|\, \lambda \in
 W(A)\}.$ The problem of characterizing linear maps on matrices or operators that preserve
numerical range or numerical radius has been studied by many
authors, see for example \cite{Chan4,Chan5, CH2, LT2} and  the
references therein. In recent years, interest in characterizing
general (non-linear) preservers of numerical ranges or numerical
radius has been growing
(\cite{BH8,BHX2,CLS7,CH1,CH3, DHHK,HwL,HHZ,HHZ2,HD,HLQ,LS7,
LPS}).

Let ${\mathcal B}_s(H)$
and  ${\mathcal B}(H)$  the space of all self-adjoint operators and
the algebra of all bounded linear operators on complex Hilbert space $H$, respectively.
 Suppose that $\mathcal A=\mathcal B(H)$ or ${\mathcal B}_s(H)$, and $F$ is the numerical range $W$ or numerical radius $w$. Let $A\circ B$
denote any  product  of  a pair of $A,B\in \mathcal A$  such as
operator product $AB$, Jordan product $AB+BA$, Jordan semi-triple
product $ABA$ and Lie product $AB-BA$. A map $\Phi: \mathcal
A\rightarrow \mathcal A$ preserves numerical range (or numerical
radius) of product $\circ$  if $F=W$ (or $F=w$) and $\Phi$ satisfies
\begin{equation} F(A\circ B)=F(\Phi(A) \circ \Phi(B)) \end{equation} for all $A, B\in
\mathcal A$.

Assume that $\Phi:{\mathcal A}\to{\mathcal A}$ satisfy Eq.(1.1). For
the case $F=W$ and $\Phi$ is surjective, it was shown in \cite{HD}
that if $A\circ B=AB$ and $\mathcal A=\mathcal B(H)$, then there
exists a unitary operator $U$ such that $\Phi(A)=\epsilon UAU^*$ for
all $A\in {\mathcal A}$, where $\epsilon\in\{-1,1\}$; if $A\circ
B=ABA$ and $\mathcal A=\mathcal B(H)$, then $\Phi$ is the multiple
of a C$^*$-isomorphism (by a cubic root of unity); if $A\circ B=AB$
and $\mathcal A=\mathcal B_s(H)$, then there exists a unitary
operator $U$ such that $\Phi(A)=\epsilon UAU^*$ for all $A\in
{\mathcal A}$, where $\epsilon\in\{-1,1\} $. For the case $F=w$,
$A\circ B=AB$, $\mathcal A=\mathcal B(H)$ and $\Phi$ is surjective,
it was proved in \cite{CH1}   that there exist a unitary or
anti-unitary operator $U$ and a unit-modular functional $f: \mathcal
A\rightarrow \mathbb C$ such that $\Phi(A)=f(A)UAU^*$ for all $A\in
{\mathcal A}$. The case of $F=w$, $A\circ B=ABA$ and $\mathcal
A=\mathcal B(H)$ was dealt with in \cite{DKLLP}. For the case when
$F=w$, $A\circ B=AB$ or $ABA$ and $\mathcal A=\mathcal B_s(H)$, the
results obtained in \cite{HHZ} reveal that there is a unitary
operator or conjugate unitary operator $U$ on $H$, a sign function
$h: {\mathcal S}(H)\rightarrow \{1,-1\}$ such that
$\Phi(T)=h(T)UTU^{*}$ for any $T\in {\mathcal B}_s(H)$. Maps
preserving numerical range of Jordan product are characterized in
\cite{DHHK,HwL,  LS7}.

Recent interest is focused on characterizing non-linear maps
preserving numerical range or numerical radius of Lie product. When
$3\leq \dim H=n<\infty$,   Li, Poon and Sze \cite{LPS} proved that a
surjective map $\Phi: \mathcal B(H)\rightarrow \mathcal B(H)$
satisfies $w(\Phi(A)\Phi(B)-\Phi(B)\Phi(A))=w(AB-BA)$ for all
$A,B\in \mathcal A$ if and only if there exists a unitary matrix $U$
such that
 $$\Phi(A)=\mu_A UA^\dag U^*+\nu_AI$$
 for all $A\in \mathcal A$, where $\mu_A,\nu_A\in {\mathbb C}$ depend on
 $A$ with $|\mu_A|=1$, $(\cdot)^\dag$ stands for one of the following four maps: $A\mapsto A,A\mapsto \bar{A}, A\mapsto A^t$ and $A\mapsto A^*$.
For arbitrary dimensional space $H$ (concluding infinite and two
dimensional cases), without assumption of surjectivity,  Hou, Li and
Qi \cite{HLQ}  gave a characterization of maps on $\mathcal
 B(H)$
 preserving numerical range of Lie product.

{\bf Theorem HLQ.} {\it Let $H,K$ be complex Hilbert spaces of
dimension $\geq 2$ and $\Phi:{\mathcal B}(H)\to{\mathcal B}(K)$ be a
map of which the range contains   all operators of rank $\leq 2$.
Then  the following statements are equivalent.}

(1) {\it $\Phi$ satisfies that $ {W}([\Phi(A),\Phi(B)])= {W}([A,B])$
for any $A,B\in{\mathcal B}(H)$. }

(2) {\it  $\dim H=\dim K$, and there exist $\varepsilon \in
\{1,-1\}$, a functional $h: \mathcal B(H) \rightarrow \Bbb C$, a
unitary operator $U\in{\mathcal B}(H,K)$, and a set ${\mathcal S}$
of operators in ${\mathcal B}(H)$, which consists of operators of
the form $aP + bI$ for an orthogonal projection $P$ on $H$ if $\dim
H \ge 3$, such that either
$$\Phi(A)=\begin{cases}
\ \varepsilon UAU^*+h(A)I &  $if$ \ A \in{\mathcal B}(H)\setminus \mathcal S,\cr
-\varepsilon UAU^* + h(A)I & $if$ \ A \in \mathcal S,\cr\end{cases}$$
or
$$\Phi(A) = \begin{cases}
\ i\varepsilon UA^tU^*+h(A)I &  $if$ \ A\in{\mathcal B}(H) \setminus
\mathcal S,\cr -i\varepsilon UA^tU^*+h(A)I &   $if$ \ A\in\mathcal S, \cr
\end{cases}$$ where $A^t$ is the transpose of $A$ with respect to an
orthonormal basis of $H$.}\\

An interesting open question is how to characterize non-linear maps
on self-adjoint operators preserving numerical radius or numerical
range of Lie product. In this paper, we solve this question for the
case when the underline space $H$ is separable.

Let $H$ be a complex separable Hilbert space and $\Phi:{\mathcal
B}_s(H)\to{\mathcal B}_s(H)$ a surjective map. Assume further that
$\dim H\geq 3$. We show that:

(a) $\Phi$ satisfies $w(AB- BA)=w(\Phi(A)\Phi(B)-\Phi(B)\Phi(A) ) $
for any $A,B\in {\mathcal B}_s(H)$ if and only if  there exist a
unitary   operator $U$ on $H$, a sign function $h: {\mathcal
B}_s(H)\rightarrow \{1,-1\}$  and a functional $f: {\mathcal
B}_s(H)\rightarrow {\Bbb R}$ such that $\Phi(T)=h(T)UTU^{*}+f(T)I$
for all $T\in {\mathcal B}_s(H)$ or $\Phi(T)=h(T)UT^tU^{*}+f(T)I$
for all $T\in {\mathcal B}_s(H)$ (See Theorem 2.1);

(b) $\Phi$ satisfies $W(AB- BA)=W(\Phi(A)\Phi(B)-\Phi(B)\Phi(A) )$
for all $A,B\in {\mathcal B}_s(H) $  if and only if there exist a
unitary operator $U$ on $H$, a scalar $\varepsilon \in\{1,-1\}$, a
subset ${\mathcal S}\subseteq {\mathcal D}(H)$, and a functional $f:
{\mathcal B}_s(H)\rightarrow {\Bbb R}$ such that $\Phi(A)=
  \varepsilon UAU^*+f(A)I$     if $A
\in{\mathcal B}_s(H)\setminus {\mathcal S}$, $\Phi(A)=-\varepsilon
UAU^* + f(A)I$ if $ A \in {\mathcal S}$, where ${\mathcal D}(H)$ is
the set of all real linear combinations of a projection and the
identity $I$ on $H$ (See Theorem 3.1).

When $\dim H=2$, not like to the maps on ${\mathcal B}(H)$ (See
Theorem HLQ list above), the maps $\Phi:{\mathcal
B}_s(H)\to{\mathcal B}_s(H)$ which preserve numerical range (radius)
of Lie product may have some other forms. Note that, in the case
$\dim H=2$ we have ${\mathcal D}(H)={\mathcal B}_2(H)$. Identifying
${\mathcal B}_s(H)$ with the space ${\bf H}_2$ of all $2\times 2$
Hermitian matrices, and define $\Psi$ on ${\bf H}_2$ by $$
\left(\begin{array}{cc} a & c+id \\ c-id &b
\end{array}\right)\mapsto\left(\begin{array}{cc} a & -c+id \\ -c-id &b
\end{array}\right).$$ It is easily checked that $\Psi$ preserves both the numerical range and the
numerical radius of Lie product. However,  no more other kind of
maps can be added in as revealed by our result. In addition, the
surjectivity assumption is not needed in the following result.

(c) A map $\Phi: {\bf H}_2\to{\bf H}_2$ preserves the numerical
radius of Lie product if and only if it preserves the numerical
range of Lie product, and in turn, if and only if there exist a
unitary matrix $U\in M_2$, a sign function $h:{\bf H}_2\to \{1,-1\}$
and a functional $f: {\mathcal B}_s(H)\rightarrow {\Bbb R}$ such
that either $\Phi(A)=
  h(A) UA^\dag U^*+f(A)I$  for all $A
\in{\bf H}_2$; or $\Phi(A)=h(A)U\Psi(A)^\dag U^* + f(A)I$ for all $
A \in {\bf H}_2$, where $(\cdot)^\dag$ is one of  the identity map
and the transpose map (See Theorem 4.1).

The paper is organized as follows. We characterize the maps preserve
the numerical radius of Lie product for the case $\dim H\geq 3$ in
Section 2 and the maps preserving numerical range of Lie product for
the case $\dim H\geq 3$ in Section 3. The last section is devoted to
the case when $\dim H=2$.

 \section{Preservers for numerical radius of Lie product}

 In the section, we devote to characterizing surjective maps on self-adjoint operators preserving numerical radius of Lie
 product for the case $\dim H\geq 3$. The following is the main result.

\begin{thm}\label{thm:1}
Let $H$ be a separable complex Hilbert space of dimension at least
three. A surjective map $\Phi\colon {\mathcal B}_s(H)\rightarrow
{\mathcal B}_s(H)$ satisfies
$$w(AB- BA)=w(\Phi(A)\Phi(B)-\Phi(B)\Phi(A) ) $$  for all $A,B\in {\mathcal
B}_s(H)$ if and only if  there exist a unitary  operator $U$ on $H$,
a sign function $h: {\mathcal B}_s(H)\rightarrow \{1,-1\}$  and a
functional $f: {\mathcal B}_s(H)\rightarrow {\Bbb R}$ such that
either
$$\Phi(T)=h(T)UTU^{*}+f(T)I$$ for all $T\in {\mathcal
B}_s(H)$; or
$$\Phi(T)=h(T)UT^tU^{*}+f(T)I$$ for all $T\in {\mathcal
B}_s(H)$. Here $T^t$ is the transpose of $T$ with respect to an
arbitrarily given orthonormal basis of $H$.\end{thm}

Before starting the proof of Theorem 2.1, we need a lemma.

{\bf Lemma 2.2.} {\it Let $H$ be a complex Hilbert space of
dimension $\geq 2$ and $A,B$ be self-adjoint operators acting on
$H$. Then the following statements are equivalent.}

(1) {\it  $w(AC-CA)=w(BC-C B)$ for every   $C\in{\mathcal B}_s(H)$.}

(2) {\it $w(AP-P A)=w(BP-P B)$ for every rank-1 projection $P$.}

(3) {\it $A+B$ or $A-B$ is a scalar.}

{\it Proof.} (3)$\Rightarrow$(1)$\Rightarrow$(2) are  obvious. Let
us check (2)$\Rightarrow$(3).

Assume (2). For any rank-1 projection $P=x\otimes x$, write
$Ax=\alpha x+\beta y$, where normalized $y$ is orthogonal to $x$.
Since $A$ is self-adjoint we have $\alpha=\langle
Ax,x\rangle\in\mathbb R$. Moreover, by self-adjointness of $A$,
$Ax\otimes x-x\otimes x A=Ax\otimes x-x\otimes (Ax)$. So relative to
decomposition $H=[x,y]\oplus H_1$, the rank-2 operator $Ax\otimes
x-x\otimes x A$ is represented by a matrix $$\left(\begin{matrix}
0 & -\bar{\beta} \\
\beta& 0
\end{matrix}\right)\oplus 0,$$
and hence $W(Ax\otimes x-x\otimes x A)=i[-|\beta|,|\beta|]$ and
$w(Ax\otimes x-x\otimes x A)=|\beta|$.
Decomposing likewise $Bx=\alpha' x+\beta' z$ we obtain $W(Bx\otimes
x-x\otimes xB) = i[-|\beta'|,|\beta'|]$ and the numerical radius of
$[B, x\otimes x]$ is $|\beta'|$. Hence by (2) we obtain
$|\beta|=|\beta'|$.

Since $A$ is self-adjoint and $\alpha=\langle Ax,x\rangle \in
\mathbb R$, it follows that $|\beta|^2=\|(Ax-\langle Ax,x\rangle\,
x)\|^2=\langle(Ax-\langle Ax,x\rangle\, x)\,,\,(Ax-\langle
Ax,x\rangle\, x)\rangle
=\langle A^2x,x\rangle- \langle Ax,x\rangle   ^2$. Similarly, for
$B$ we obtain $|\beta'|^2  =\langle B^2x,x\rangle-  \langle
Bx,x\rangle  ^2$. If follows from $|\beta|^2=|\beta'|^2$   that
\begin{equation}\label{betabeta}
\langle A^2x,x\rangle-\langle B^2x,x\rangle= \langle Ax,x\rangle
 ^2- \langle Bx,x\rangle  ^2
\end{equation}
for every normalized vector $x$. Let $y,z$ be two orthogonal
normalized vectors. Then $x=\frac{\sqrt{2}}{2} (e^{i\xi}y+z) $ is
also normalized  for every $\xi\in[-\pi,\pi]$. After inserting $x$
in Eq.(2.1)  we obtain
\begin{equation}\label{Fourier1}
\begin{aligned}
 0  =&2\langle A^2(e^{i\xi}y + z),e^{i\xi}y + z\rangle - 2\langle B^2(e^{i\xi}y + z),e^{i\xi}y + z\rangle \\
 & -  \bigl(\langle A(e^{i\xi}y + z), e^{i\xi}y + z\rangle
\bigr)^2+\bigl(\langle B(e^{i\xi}y + z),e^{i\xi}y + z\rangle
\bigr)^2.
\end{aligned}
\end{equation}
Taking only the coefficient at $e^{2 i\xi}$ in the expansion of
Eq.(2.2) in a Fourier series,   Eq.(2.2) reduces into
\begin{equation}\label{eq:langle..rangle^2}
\langle By,z\rangle^2=\langle A y,z\rangle^2
\end{equation}
for every pair of orthonormal $y,z$.  So, for any $x\in H$ and
$f\in[Ax,x]^\perp$, we have $\langle Bx,f\rangle=0$. This entails
that $Bx\in[Ax,x]$.  Thus, for any $x\in H$, there exist
$\alpha_x,\beta_x\in\mathbb{C}$ such that $Bx=\alpha_xAx+\beta_xx$.
By Eq.(2.3), we have $\langle
Ax,f\rangle^2=\langle\alpha_xAx,f\rangle^2 =\alpha_x^2\langle
Ax,f\rangle^2$ holds for all $f\in[x]^\perp$, which implies that
  $\alpha_x=\pm1$. It follows from
$|\beta_x|\|x\|\leq\|Bx\|+\|\alpha_xAx\|\leq(\|B\|+\|A\|)\|x\|$ that
$|\beta_x|\leq\|B\|+\|A\|$. Therefore, $B$ is a regular local linear
combination of $A$ and $I$, and then, by \cite{Hou1}, $B$ is a
linear combination of $A$ and $I$. So $B=\alpha A+\beta I$ with
$\alpha\in\{-1,1\}$ and $\beta\in{\mathbb R}$, as desired.
 \hfill$\Box$

\begin{proof}[Proof of Theorem 2.1] The ``if'' part is obvious, we check the ``only if'' part.

Assume first that $\Phi$ is injective; then, $\Phi$ is bijective.
 Clearly $\Phi$ preserves zeros of Lie product.  So, by \cite{Mol-Sem},  there exists a unitary or conjugate unitary operator $U$ such that, for any
 rank-1 positive operator  $P=x\otimes x$ with  unit  vector $x\in H$, we have
 $$\Phi(P)=U(\lambda_P P +\mu_P I)U^* $$
for some $ \lambda_P, \mu_P \in \Bbb R$. Without loss of generality
we can assume in the sequel that $U=I$.
 Taking any unit vectors $x,y$, which are orthogonal to each other,
and let  $Q=y\otimes y$ and $Z=(x+y)\otimes (x+y)$. It easily
follows that in the orthogonal decomposition  $H=[x,y]\oplus H_1$,
where $[x,y]$ stands for the subspace spanned by $\{x,y\}$, we have
$PZ-ZP= \left(\begin{matrix}
0 & 1\\
-1 &0
\end{matrix}\right)\oplus 0$, whose numerical range is $[-i,i]$ and so numerical radius is 1. The same conclusion
holds for the numerical range of $QZ-ZQ$. Comparing the numerical radius of  $PZ-ZP$ and $\Phi(P)\Phi(Z)-\Phi(Z)\Phi(P)$, we obtain
$$1=w(\Phi(P)\Phi(Z)-\Phi(Z)\Phi(P))=|\lambda_P\lambda_Z| w(PZ-ZP)=|\lambda_P\lambda_Z|.$$ This is possible only if $\lambda_P\lambda_Z=\pm1$
since $\lambda_P, \lambda_Z\in \mathbb R$. Similarly we have
$\lambda_Q\lambda_Z=\pm 1$. Hence $\lambda_P=\pm\lambda_Q$ for
orthogonal $P, Q$. Now, given any rank-one self-adjoint operator
$R$, there exists a rank-one self-adjoint operator $T$ which is
orthogonal to $R$ and $P$. Similar to the above discussion, we have
$\lambda_T=\pm \lambda_R$ and $\lambda_T=\pm \lambda_P$, so
$\lambda_P=\pm\lambda_R$ for any $P, R$. It follows that
$\lambda_P=\pm 1$. \if Write $\lambda_P$ as $\lambda$. And as
$\lambda_P\lambda_Z=\pm 1$, we deduce $\lambda=\pm 1$ for any
rank-one self-adjoint operator $P$.\fi


\if Now for rank-one projection $P=x\otimes x$,  assume
$\Phi(P)=\lambda_P P +\mu_P I$.\fi Now, for arbitrary self-adjoint
$A$,
\begin{equation}\label{eq:rank-one}
w(Ax\otimes x-x\otimes x A)=|\lambda_P| w(\Phi(A)x\otimes x-x\otimes
x\Phi(A))=w(\Phi(A)x\otimes x-x\otimes x\Phi(A))
\end{equation}
holds for every rank-1 projection $P=x\otimes x$. \if Decompose
$Ax=\alpha x+\beta y$, where normalized $y$ is orthogonal to $x$.
Since $A$ is self-adjoint we have $\alpha=\langle
Ax,x\rangle\in\mathbb R$. Moreover, by self-adjointness of $A$,
$Ax\otimes x-x\otimes x A=Ax\otimes x-x\otimes (Ax)$. So relative to
decomposition $H=[x,y]\oplus H_1$, the rank-two operator $Ax\otimes
x-x\otimes x A$ is represented by a matrix
$$\left(\begin{matrix}
0 & -\bar{\beta} \\
\beta& 0
\end{matrix}\right)\oplus 0,$$
and hence $W(Ax\otimes x-x\otimes x A)=i[-|\beta|,|\beta|]$ and $w(Ax\otimes x-x\otimes x A)=|\beta|$.
Decomposing likewise $\Phi(A)x=\alpha' x+\beta' z$ we obtain $W(\Phi(A)x\otimes x-x\otimes x\Phi(A)) =
i[-|\beta'|,|\beta'|]$ and the numerical radius of $[\Phi(A), x\otimes x]$ is $|\beta'|$. Hence
by
Eq.(2.1) we obtain  $|\beta|=|\beta'|$.

Since $A$ is self-adjoint and $\alpha=\langle Ax,x\rangle \in
\mathbb R$, it follows that $|\beta|^2=\|(Ax-\langle Ax,x\rangle\,
x)\|^2=\langle(Ax-\langle Ax,x\rangle\, x)\,,\,(Ax-\langle
Ax,x\rangle\, x)\rangle
=\langle A^2x,x\rangle-\bigl(\langle Ax,x\rangle \bigr)^2$.
Similarly, for $\Phi(A)$ we obtain $|\beta'|^2  =\langle
\Phi(A)^2x,x\rangle-\bigl(\langle \Phi(A)x,x\rangle \bigr)^2$. If
follows from $|\beta|=|\beta'|$   that
\begin{equation}\label{betabeta}
\langle A^2x,x\rangle-\langle \Phi(A)^2x,x\rangle=\bigl(\langle Ax,x\rangle
\bigr)^2-\bigl(\langle \Phi(A)x,x\rangle \bigr)^2
\end{equation}
for every normalized vector $x$.
Let $y,z$ be two orthogonal normalized vectors. Then $x=\frac{\sqrt{2}}{2} (e^{i\xi}y+z) $ is also normalized  for every
$\xi\in[-\pi,\pi]$. After inserting $x$ in Eq.(2.2)  we obtain
\begin{equation}\label{Fourier1}
\begin{aligned}
 0  =&2\langle A^2(e^{i\xi}y + z),e^{i\xi}y + z\rangle - 2\langle \Phi(A)^2(e^{i\xi}y + z),e^{i\xi}y + z\rangle \\
 & -  \bigl(\langle A(e^{i\xi}y + z), e^{i\xi}y + z\rangle
\bigr)^2+\bigl(\langle \Phi(A)(e^{i\xi}y + z),e^{i\xi}y + z\rangle \bigr)^2.
\end{aligned}
\end{equation}
Taking only the coefficient at $e^{2 i\xi}$ in the expansion of Eq.(2.3) in a Fourier series, identity Eq.(2.3) reduces into
\begin{equation}\label{eq:langle..rangle^2}
\langle \Phi(A)y,z\rangle^2=\langle A y,z\rangle^2
\end{equation}
for every pair of orthonormal $y,z$. From this we infer that also
$|\langle \Phi(A)y,z\rangle|=|\langle A y,z\rangle|$ wherefrom, by
\cite[Theorem 2.2]{kllpr},\fi By Lemma 2.2, $\Phi(A)=\lambda_A
A+\delta_A I$ for some scalar $\lambda_A\in\{-1,1\}$ and some scalar
$\delta_A$.

Finally we show that one only needs the surjective assumption. Here
we borrow an idea from \cite{LPS}. If $\Phi(A)=\Phi(B)$, then
$$\begin{array}{rl} w(AC-CA)=& w(\Phi(A)\Phi(C)-\Phi(C)\Phi(A))\\ =& w(\Phi(B)\Phi(C)-\Phi(C)\Phi(B))=w(BC-CB) \end{array}$$
for all $C\in{\mathcal B}_s(H)$. By Lemma 2.2 we get $B=\alpha
A+\beta I$ for some $\alpha\in\{-1,1\}$ and $\beta\in{\mathbb R}$.
On the other hand, for any $A$, there is some $D$ such that
$\Phi(D)=-\Phi(A)$, which gives
$w(DC-CD)=w(\Phi(D)\Phi(C)-\Phi(C)\Phi(D))=w(\Phi(A)\Phi(C)-\Phi(C)\Phi(A))=w(AC-CA)$
for all $C$. Again by Lemma 2.2, we get $D=\lambda A+\gamma I$ for
some $\lambda\in\{-1,1\}$ and $\gamma\in{\mathbb R}$. For any
$A,B\in{\mathcal B}_s(H)$, we say $A\sim B$ if $w(AC-CA)=w(BC-CB)$
for all $C\in{\mathcal B}_s(H)$. By Lemma 2.2, $\sim$ is an
equivalent relation and $A\sim B$ if and only if $B=\alpha A+\beta
I$ for some $\alpha\in\{-1,1\}$ and $\beta\in{\mathbb R}$. Let
${\mathcal E}_A=\{B\in{\mathcal B}_s(H): B\sim A\}$. For each
equivalent class ${\mathcal E}_A$ pick a representative, for example
$A$, and write ${\mathcal A}$ the set of these representatives.
Since $\Phi$ is surjective, for each $A\in{\mathcal A}$, ${\mathcal
E}_A$ and $\Phi^{-1}({\mathcal E}_A)$ have the same cardinality $c$.
Thus  there exists a map $\Psi:{\mathcal B}_s(H)\to {\mathcal
B}_s(H)$ which maps bijectively $\Phi^{-1}({\mathcal E}_A)$ onto
${\mathcal E}_A$ for each  $A\in{\mathcal A}$. Obviously, $\Psi$ is
bijective and $\Psi(A)\sim\Phi(A)$ for all $A\in{\mathcal B}_s(H)$.
Then
$$w(\Psi(A)\Psi(B)-\Psi(B)\Psi(A))=w(\Phi(A)\Phi(B)-\Phi(B)\Phi(A))=w(AB-BA)$$
for all $A,B\in{\mathcal B}_s(H)$. By the previous part of our proof
of the theorem under the bijective assumption, $\Psi$ has the
desired form, and hence $\Phi$ has the desired form as $\Phi(A)\sim
\Psi(A)$. So Theorem 2.1 holds true,
 completing the proof. \end{proof}

\section{Preservers for numerical range of Lie product }

This section is devoted to characterizing maps that preserve the
numerical range of Lie product of self-adjoint operators. Our main
result is Theorem 3.1, which is not a direct corollary of Theorem
2.1 for numerical radius preservers, since much more effort should
be paid to determine the structure of the sign function $h:
{\mathcal B}_s(H)\rightarrow \{1,-1\}$.

Denote $\mathcal D$ the set of all real linear combinations of a
projection and the identity $I$, that is, $ {\mathcal D} = \{ \alpha
P + \delta I: P \mbox{ is a projection in} \ {\mathcal B}_s(H),
\alpha, \delta \in \mathbb R \} \subset {\mathcal B}_s(H)$. It is
clear that $\mathcal D$ is those self-adjoint operators that are
also quadric algebraical operators.

\begin{thm}\label{thm:2}
Let $H$ be a complex separable Hilbert space of dimension at least
3. A surjection $\Phi\colon {\mathcal B}_s(H)\rightarrow {\mathcal
B}_s(H)$ satisfies
$$W(AB- BA)=W(\Phi(A)\Phi(B)-\Phi(B)\Phi(A) )$$
for all $A,B\in {\mathcal B}_s(H) $  if and only if there exist a
unitary operator $U$ on $H$,  a scalar $\varepsilon \in\{1,-1\}$, a
set ${\mathcal S}\subseteq \mathcal D$, and a functional $f:
{\mathcal B}_s(H)\rightarrow {\Bbb R}$ such that
$$\Phi(A)=\left\{ \begin{array}{lll}
\ \varepsilon UAU^*+f(A)I &  {\rm if} & A \in{\mathcal
B}_s(H)\setminus {\mathcal S},\\ -\varepsilon UAU^* + f(A)I & {\rm
if} & A \in {\mathcal S}.\end{array}\right.$$

\end{thm}

To prove the above result we need a lemma, which gives a
characterization of the quadric algebraic self-adjoint operators,
that is, the operators in $\mathcal D$, inn terms of the numerical
range of Lie product.

 {\bf Lemma 3.2.} {\it Let $H$ be a complex Hilbert space
with $\dim H\geq 3$ and $A\in{\mathcal B}_s(H)$. Then the following
statements are equivalent.}

(1) $A\in{\mathcal D}$.

(2) {\it  $W(AB-BA)=-W(AB-BA)$ for all $B\in{\mathcal B}_s(H)$.}

(3) {\it $W(AB-BA)=-W(AB-BA)$ for all $B\in{\mathcal B}_s(H)$ of
rank $\leq 2$.}

{\it Proof.}  (1)$\Rightarrow$(2). Assume $A\in{\mathcal D}$, then
$A=\alpha P+\gamma I$ for some projection $P$ and some scalars
$\alpha, \gamma\in{\mathbb R}$.  As the case $A=\alpha I$ is
obvious, we may assume that, there exists a space
decomposition $H=H_1\oplus H_2$ such that $A=\left(\begin{array}{cc} \alpha I_{H_1} & 0\\
0 & \beta I_{H_2}
\end{array}\right)$ with $\dim H_i>0$, $i=1,2$, and
$\alpha\not=\beta$. For any $B=\left(\begin{array}{cc} B_{11} & B_{12} \\
B_{12}^* & B_{22} \end{array}\right)\in {\mathcal B}_s(H_1\oplus
H_2)$,
$AB-BA=(\alpha-\beta)\left(\begin{array}{cc} 0 & B_{12} \\
-B_{12}^* & 0 \end{array}\right)$. Let $U=\left(\begin{array}{cc} I_{H_1} & 0 \\
0 & -I_{H_2} \end{array}\right)$; then $U$ is unitary and
$U(AB-BA)U^*=-(AB-BA)$. So, we always have $W(AB-BA)=-W(AB-BA)$,
that is, (2) is true.

(2)$\Rightarrow$(3) is obvious.

(3)$\Rightarrow$(1). Note that $A\in{\mathcal B}_s(H) \setminus
\mathcal D$ if and only if the spectrum $\sigma(A)$ has at least
three points, and in turn, if and only if there exists a vector $x$
such that $\{x, Ax, A^2 x\}$ is linearly independent.  For such $x$,
take an orthonormal basis $\{e_1,e_2,e_3\}$ of $[x,Ax,A^2x]$ with
$e_1\in [x]$ and $e_2\in[x,Ax]$. Then, with respect to the space
decomposition $H=[e_1]\oplus [e_2]\oplus [e_3]\oplus
\{e_1,e_2,e_3\}^\perp$, $A$ has the matrix representation of the
form
$$A=\left(\begin{array}{cccc} a_{11} & a_{21} & 0 & 0 \\
a_{21} & a_{22} & a_{32} & 0 \\ 0 & a_{32} & a_{33} & A_{34}\\
0 & 0 & A_{34}^* &A_{44} \end{array}\right) $$ with $a_{11},
a_{22},a_{33}$ real numbers, $a_{21}>0$, $a_{32}>0$ and
$A_{44}=A_{44}^*$. Let
$$B=\left(\begin{array}{cccc} 1 & \beta & 0 & 0 \\
\bar{\beta} &0 & 0 & 0 \\ 0 & 0 &0 & 0\\
0 & 0 & 0 &0 \end{array}\right) $$ with ${\rm
Im}\beta=\frac{1}{2i}(\beta-\bar{\beta})\not=0$. Then, $B$ is of
rank two and
$$AB-BA=\left(\begin{array}{cccc} -2({\rm Im}\beta)a_{21} & -a_{21}+\beta(a_{11}-a_{22}) & -\beta a_{32} & 0 \\
a_{21}-\bar{\beta}(a_{11}-a_{22}) & 2({\rm Im}\beta)a_{21} & 0 & 0 \\ \bar{\beta}a_{32} & 0 & 0& 0\\
0 & 0 & 0 &0 \end{array}\right), $$ which is a rank-3 skew
self-adjoint operator with zero trace. Clearly,
$W(AB-BA)\not=-W(AB-BA)$. Hence (3) implies (1). \hfill$\Box$

\if {\bf Lemma 2.5.} {\it Let $H$ be a complex Hilbert space with
$\dim H\geq 3$ and $A\in{\mathcal B}_s(H)$. Then the following
statements are equivalent.}

(1) $A=\alpha P+\gamma I$ for some rank-1 projection $P$ and
$\alpha, \gamma\in{\mathbb R}$.

(2) {\it  $W(AB-BA)=-W(AB-BA)$ is symmetric relative to 0 for all
$B\in{\mathcal B}_s(H)$.}

(3) {\it $W(AB-BA)=-W(AB-BA)$ is symmetric relative to 0 for all
$B\in{\mathcal B}_s(H)$ of rank $\leq 3$.}

{\it Proof.} (1)$\Rightarrow$(2). If $A=\alpha x\otimes x+\gamma I$,
then, for any $B\in{\mathcal B}_s(H)$, $AB-BA$ is a skew
self-adjoint operator of rank $\leq 2$. As ${\rm Tr}(AB-BA)=0$, we
have $\sigma(AB-BA)=\{0,-ai,ai\}$ for some nonnegative real number
$a$. Thus $W(AB-BA)=[-ai,ai]$ is symmetric relative to the origin 0.

(2)$\Rightarrow$(3) is clear.

(3)$\Rightarrow$(1). By Lemma 2.4, $A\in{\mathcal P}$, that is,
$A=\alpha P+\gamma I$ for some projection $P$ and real numbers
$\alpha, \gamma$. We have to show that $P=x\otimes x$ for some unit
vector $x$.

Without loss of generality, we may assume that $\alpha\not=0$ and
$P\not=0,I$. Assume, on the contrary, that $A\not\in{\mathbb
C}{\mathcal P}_1+{\mathbb C}I$; then $\dim H\geq 4$ and there exists
space decomposition $H=H_1\oplus H_2\oplus H_3$ such that
$$A=A_1\oplus \alpha I_{H_2}\oplus \beta I_{H_3},$$
where $\alpha\not=\beta$ and
$$A_1=\left(\begin{array}{cccc} \alpha & 0 & 0 & 0\\ 0 & \alpha & 0
&0 \\ 0 &0 &\beta & 0\\ 0&0&0&\beta
\end{array}\right)=\left(\begin{array}{cc} \alpha I_2 &0_2\\ 0_2 & \beta I_2\end{array}\right).$$
Let $B=B_1\oplus 0\oplus 0$ with \if $B_1=\left(\begin{array}{cc}  B_{11}&B_{12}\\
B_{21} & B_{22}\end{array}\right)\in{\mathcal B}(H_1)$. Then
$$[A,B]=[A_1, B_1]\oplus 0\oplus
0=(\alpha-\beta)\left(\begin{array}{cc} 0_2 &B_{12}
\\-B_{21}&0_2\end{array}\right)\oplus 0\oplus 0.$$ Write $B_{12}=\left(\begin{array}{cc}  \alpha_1&\alpha_2\\
\alpha_3 &\alpha_4\end{array}\right)$ and $-B_{21}=\left(\begin{array}{cc}  \beta_1&\beta_2\\
\beta_3 &\beta_4\end{array}\right)$. It is easily checked that
$[A_1, B_1]$ has four eigenvalues
$$\begin{array}{rl} \lambda=&\pm (\alpha-\beta)\frac{1}{\sqrt{2}}  [\alpha_1\beta_1-\alpha_4\beta_4+\alpha_2\beta_3-\alpha_3\beta_2 \\
 & \pm\sqrt{(\alpha_1\beta_1-\alpha_4\beta_4+\alpha_2\beta_3-\alpha_3\beta_2)^2
-4(\alpha_1\alpha_4-\alpha_2\alpha_3)(\beta_1\beta_4-\beta_2\beta_3)}]^{\frac{1}{2}}.\end{array}
$$
If we take $\alpha_i$ and $\beta_i$ so that
$\alpha_1\beta_1-\alpha_4\beta_4+\alpha_2\beta_3-\alpha_3\beta_2=0$
and
$(\alpha_1\alpha_4-\alpha_2\alpha_3)(\beta_1\beta_4-\beta_2\beta_3)=-d<0$,
then $$\sigma((\alpha-\beta)\left(\begin{array}{cc} 0_2 &B_{12}
\\-B_{21}&0_2\end{array}\right))=\{\pm (\alpha-\beta)\sqrt{d}, \pm
i(\alpha-\beta)\sqrt{d}\}.$$ Thus we have
$$\sigma([A,B])=\{0, \pm (\alpha-\beta)\sqrt{d}, \pm
i(\alpha-\beta)\sqrt{d}\}.$$  Now, take \fi $$B_1=(\alpha-\beta)^{-1}\left(\begin{array}{cccc} 2 & 1 & 1 & 2\\ -1 &-1 &-1 &-1 \\
2 &1 &1 &2\\ -1&-1&-1&-1
\end{array}\right);$$ we see that  rank$(B )=2$,   $$ [A_1,B_1]=\left(\begin{array}{cccc} 0 & 0 & 1 & 2\\ 0 &0 &-1 &-1 \\
-2 &-1 &0 &0\\ 1&1&0&0
\end{array}\right) $$ and  $\sigma([A_1,B_1])=\{\pm 1,
\pm i\}$. Note that $W([A,B])=W([A_1,B_1])$ as $0\in W([A_1,B_1])$.
Thus, if $W([A,B])$ is an elliptic disc centered at 0, then by Lemma
2.2, $W([A,B])$ is either an elliptic disc of foci $\{-1,1\}$, or an
elliptic disc of foci $\{-i,i\}$. This forces that $W([A,B])$ is a
circular disc centered at 0 because  $\sigma ({\rm
Re}([A_1,B_1])=\sigma ({\rm
Im}([A_1,B_1])=\{-1.6283,-0.9212,0.9212,1.6283\}$ implies the major
axis and the minor axis have the same length.   However,
$\sigma({\rm Re}
(e^{\frac{i\pi}{3}}[A_1,B_1]))=\{-1.5939,-0.9795,0.9795,1.5939\}$,
which means that $W({\rm Re}
(e^{\frac{i\pi}{3}}[A_1,B_1]))\not=W({\rm Re}([A_1,B_1])$ and then
$W([A,B])$ cannot be   a circular disc, a contradiction. Hence the
numerical range of $[A,B]$ cannot be an elliptic disc with center 0
and $A$ does not meet the condition (3). \hfill$\Box$ \fi

{\bf Proof of Theorem 3.1.}

Assume $\dim H\geq 3$. Then $\Phi$ satisfies the assumption of
Theorem 2.1, and hence  there exist a unitary operator or conjugate
unitary operator $U$ on $H$, a sign function $h: {\mathcal
B}_s(H)\rightarrow \{1,-1\}$ and a functional $f: {\mathcal
B}_s(H)\rightarrow {\Bbb R}$ such that $\Phi(T)=h(T)UTU^{*}+f(T)I$
for all $T\in {\mathcal B}_s(H)$.

We assert that the case   $U$ is a conjugate unitary operator cannot
occur. Assume on the contrary that $\Phi(T)=h(T)UTU^{*}+f(T)I$ for
any $T\in {\mathcal B}_s(H)$, where $U$ is conjugate. Take
arbitrarily an orthonormal basis of $H$, one sees that there exists
a unitary operator $V$ such that $\Phi(T)=h(T)VT^tV^{*}+f(T)I$ for
any $T\in {\mathcal B}_s(H)$, where $T^t$ is the transpose of $T$
with respect to the given basis. Thus we have
\begin{equation}\label{betabeta}
\begin{aligned}
W(AB-BA)&=W(\Phi(A)\Phi(B)-\Phi(B)\Phi(A)) \\
&=h(A)h(B)W(VA^tB^tV^*-VB^tA^tV^*)\\ &=h(A)h(B)W((BA-AB)^t)\\
&=-h(A)h(B)W(AB-BA).
\end{aligned}
\end{equation}

Let $\{x,y,z\}$ be an orthonormal set of $H$ and consider the space
decomposition $H=[x,y,z]\oplus [x,y,z]^\perp$. For any scalars
$\alpha,\beta,\gamma$ with
$\alpha\beta\bar{\gamma}-\bar{\alpha}\bar{\beta}\gamma\not=0$, and
any real numbers $b_{11},b_{22}, b_{33}$, let
\begin{equation}\label{betabeta}
\begin{aligned}  B=\left(\begin{array}{ccc} b_{11} &\alpha &\gamma \\
\bar{\alpha} & b_{22} & \beta \\
\bar{\gamma} & \bar{\beta} & b_{33} \end{array}\right)\oplus 0\in
{\mathcal B}_s(H).
\end{aligned}
\end{equation}
 Then for any self-adjoint operator of the form
\begin{equation}\label{betabeta}
\begin{aligned} A=\left(\begin{array}{ccc} a_{1} &0 &0 \\
0 & a_{2} & 0 \\
0 & 0 & a_{3} \end{array}\right)\oplus A_2  \end{aligned}
\end{equation}
 with distinct
$a_1,a_2,a_3$, we have $AB-BA=C_1\oplus 0$, where
$$C_1=\left(\begin{array}{ccc} 0 &(a_1-a_2)\alpha &(a_1-a_3)\gamma \\
(a_2-a_1)\bar{\alpha} & 0 & (a_2-a_3)\beta \\
(a_3-a_1)\bar{\gamma} & (a_3-a_2)\bar{\beta} & 0
\end{array}\right).$$
As $\det
(C_1)=(a_1-a_2)(a_2-a_3)(a_3-a_1)(\alpha\beta\bar{\gamma}-\bar{\alpha}\bar{\beta}\gamma)\not=0$,
$\sigma(C_1)=\{it_1,it_2,it_3\}$ with $t_i\not=0$, $i=1,2,3$,
$t_1\leq t_2\leq t_3$ and $t_1+t_2+t_3=0$. So $W(AB-BA)=i[t_1,t_3]$
and $t_1\not=-t_3$. By Eq.(3.1) we obtain that
$$-h(A)h(B)[it_1,it_3]=[it_1,it_3]$$
and this forces $-h(A)h(B)=1$. If $h(B)=-1$, then $h(A)=1$ for all
$A$ of the    form in Eq.(3.3) and $h(B)=-1$ for all $B$ of the form
in Eq.(3.2). Consequently, $h(B)h(B')=1$ for any $B,B'$ of the form
in Eq.(3.2).

Now take self-adjoint operators of rank two
$$B=\left(\begin{array}{ccc} 0 & i & 1 \\ -i & 0 & 2 \\ 1 & 2 &
0 \end{array}\right)\oplus 0 \quad {\rm and}\quad
B'=\left(\begin{array}{ccc} 0 &1+ i & 1 \\ 1-i & 0 & 2i \\ 1 & -2i &
0
\end{array}\right)\oplus 0.$$ It is clear that $i(BB'-B'B)$ is a
rank-3 self-adjoint operator and hence
$W(BB'-B'B)\not=-W(BB'-B'B)=-h(B)h(B')W(BB'-B'B)$, contradicting to
Eq.(3.1).

So,
$$\Phi(A)=h(A)UAU^*+f(A) I$$ for all $A\in{\mathcal B}_s(H)$.

It is clear by Lemma 3.2 that $h(A)$ can take any value of $-1$ and
1 if $A\in{\mathcal D}$. So, to complete the proof, we have to show
that $h \colon {\mathcal B}_s(H) \to \{-1, 1\}$ is constant on
${\mathcal B}_s(H) \setminus \mathcal D$.

By Lemma 3.2, for any $A\in{\mathcal B}_s(H)\setminus {\mathcal D}$,
there exists rank-2 $B\in{\mathcal B}_s(H)\setminus {\mathcal D}$
such that $W(AB-BA)\not=-W(AB-BA)$. So we need only to show that
$h(A)=h(B)$ holds for any rank-two $A,B\in{\mathcal B}_s(H)\setminus
{\mathcal D}$.

{\bf Claim 1.} For any orthonormal set $\{x,y,z\}$ and any nonzero
real numbers $a,b,c,d,e,f$ with $a\not=b,c\not=d$ and $e\not=f$, we
have $h(ax\otimes x+by\otimes y)=h(cx\otimes x+dz\otimes
z)=h(ey\otimes y+fz\otimes z)$.

Assume $A$ is a rank-2 self-adjoint not in ${\mathcal D}$. Then
there exist orthonormal $x,y\in H$ and nonzero distinct real numbers
$a,b$ such that $A=\left(\begin{array}{cc} a &0 \\ 0 & b
\end{array}\right)\oplus 0$ with respect to the space decomposition
$H=[x,y]\oplus [x,y]^\perp$. Take arbitrarily two unit vectors
$z,z'\in[x,y]^\perp$ and nonzero complex numbers $\alpha,\beta,
\gamma, \alpha',\beta', \gamma'$ so that ${\rm
Re}(\alpha\beta\bar{\gamma}) =0$ and ${\rm
Re}(\alpha'\beta'\bar{\gamma}') =0$, and let
$B=B(x,y,z;\alpha,\beta,\gamma)={\rm Re} (x\otimes (\alpha y+\gamma
z)+\beta y\otimes z)$, $B'=B(x,y,z'; \alpha',\beta',\gamma')={\rm
Re} (x\otimes (\alpha' y+\gamma' z')+\beta' y\otimes z')$. Then $A$
has the form in Eq.(3.1) and $B,B'$ have the form in Eq.(3.2). By
what proved previously, we see that both $B,B'$ are of rank-2 and
$h(B)=h(A)=h(B')$ as $W(AB-BA)\not=-W(AB-BA)$ and
$W(AB'-B'A)\not=-W(AB'-B'A)$. It is also clear that
$$h(ax\otimes x+by\otimes y)=h(B(x,y,z;\alpha ,\beta ,\gamma ))=h(B(
\pi(x,y,z) ;\alpha_1,\beta_1,\gamma_1))$$ holds for any permutation
$\pi(x,y,z)$ of $(x,y,z)$ and any nonzero numbers
$\alpha_1,\beta_1,\gamma_1$ with ${\rm
Re}\alpha_1\beta_1\bar{\gamma}_1=0$. For example,
$$h(ax\otimes x+by\otimes
y)=h(B(x,y,z;\alpha,\beta,\gamma))=h(B(z,x,y;\alpha_1,\beta_1,\gamma_1)).$$
It follows that \begin{equation}h(ax\otimes x+by\otimes
y)=h(cx\otimes x+dz\otimes z)=h(ey\otimes y+fz\otimes
z)\end{equation} hold for any orthonormal set $\{x,y,z\}$ and any
nonzero real numbers $a,b,c,d,e,f$ with $a\not=b,c\not=d$ and
$e\not=f$. So Claim 1 is true.

{\bf Claim 2.} If $\dim H\geq 4$, then $h(A)=h(B)$ holds for any
rank-2 $A,B\in{\mathcal B}_s(H)\setminus{\mathcal D}$.

Let  $A=ax\otimes x+by\otimes y$ and $B=cu\otimes u+dv\otimes v$ be
any two rank-2 self-adjoint operators that are not in $\mathcal D$,
where $x\perp y$ and $u\perp v$. $\dim H\geq 4$ implies that
$[x,y,u]\not=H$. Take $y'\in[x,y,u]^\perp$. By Claim 1 we have
$$\Phi(ax\otimes x+by\otimes y)=\Phi(by\otimes y+cy'\otimes y').$$
So, replacing $y$ by $y'$ if necessary,  we may assume that $y\perp
u$ in the sequel.

If $[x,y,u,v]\not=H$, one can pick  a unit vector
$z\in[x,y,u,v]^\perp$. Then, by Claim 1 or Eq.(3.4),
$$\begin{array}{rl} h(A)=& h(ax\otimes x+by\otimes y)=h(ay\otimes y+bz\otimes z)\\ =& h(cu\otimes u+bz\otimes
z)=h(cu\otimes u+dv\otimes v)=h(B).\end{array}$$

If $ [x,y,u,v]= H$, then $\dim H= 4$. Take unit vectors $z \in[x,y,u
]^\perp$ and $z'\in[y,u,v]^\perp$. Applying Claim 1 again, we see
that
$$\begin{array}{rl}h(A)=& h(ax\otimes x+by\otimes y)=h(ay\otimes
y+bz\otimes z)\\= &h(ay\otimes y+bu\otimes u) =h(cu\otimes u
+dz'\otimes z')\\ = & h(cu\otimes u +dv\otimes v)=h(B).\end{array}$$

Finally, let us consider the case $\dim H=3$.

{\bf Claim 3.} If $\dim H=3$, then $h(A)=h(B)$ holds for any rank-2
$A,B\in{\mathcal B}_s(H)\setminus{\mathcal D}$.

 Assume that $\dim H=3$ and write $A=ax\otimes x+by\otimes y$ and $B=cu\otimes
u+dv\otimes v$, where $x\perp y, u\perp v$. If $[x,y,u,v]\not=H$,
then $[x,y]=[u,v]$. It is obvious $h(A)=h(B)$ whenever $u$ is
linearly dependent to $x$ or $y$. So we may assume that
$u,v\not\in[x]\cup [y]$.  Pick a unit vector $z\in [x,y]^\perp$. By
Claim 1 we see that $h(A)=h(ax\otimes x+by\otimes y)=h(az\otimes
z+by\otimes y)$ and $h(B)=h(cu\times u+dv\otimes v)=h(cz\otimes
z+dv\otimes v)$. It reduces to consider $A'= az\otimes z+by\otimes
y$, $B'=cz\otimes z+dv\otimes v$. Note that $[z,y,v]=H$. So we may
always require that $ [x,y,u,v]= H$.

 Take unit vector
$z\in[x,y]^\perp$; then $A$ and $B$ have matrix representations
$$A=\left(\begin{array}{ccc} a &0 &0 \\ 0 &b & 0\\
0&0&0\end{array}\right)\quad\mbox{\rm and}\quad B=\left(\begin{array}{ccc} \xi_{1} &\alpha &\gamma \\ \bar{\alpha} &\xi_{2}& \beta\\
\bar{\gamma}&\bar{\beta}&\xi_{3}\end{array}\right)$$ with $a,b,0$
are distinct to each other, $B$ has three distinct eigvalues,
$(\gamma,\beta,\xi_{3})\not=(0,0,0)$. If $\alpha=\beta=\gamma=0$ or
${\rm Im} (\alpha\beta\bar{\gamma})\not=0$, clearly we already have
$h(B)=h(A)$ (see the argument after Eqs. (3.3)-(3.4)).

In the sequel assume that $(\alpha,\beta,\gamma)\not=(0,0,0)$ but
${\rm Im} (\alpha\beta\bar{\gamma})=0$.

{\bf Subcase 1.} Two of $\alpha,\beta,\gamma$ are 0.

Without loss of generality, say $\beta=\gamma=0$. Then
$$B=\left(\begin{array}{ccc} \frac{\bar{\alpha}}{k} &\alpha &0 \\ \bar{\alpha} &k\alpha & 0\\
0&0&\xi_3\end{array}\right)$$ for some $k\not=0$ as rank$B=2$ and
$\xi_3\not=0$. Let
$$C_{t,s}=\left(\begin{array}{ccc} 0 &t &i \\ t &0 & s\\
-i&s&0\end{array}\right)$$ for nonzero $t,s\in{\mathbb R}$. By the
previous discussion we have $h(A)=h(C_{t,s})$. Consider
$$BC_{t,s}-C_{t,s}B=\left(\begin{array}{ccc} t(\alpha-\bar{\alpha}) &t(\frac{\bar{\alpha}}{k}-k\alpha) &\frac{i\bar{\alpha}}{k}+s\alpha-i\xi_3 \\
 -t(\frac{\bar{\alpha}}{k}-k\alpha) &- t(\alpha-\bar{\alpha})  & i\bar{\alpha}+sk\alpha-s\xi_3\\
\frac{i\bar{\alpha}}{k}-s\bar{\alpha}-i\xi_3 &
i\alpha-sk{\alpha}+s\xi_3&0\end{array}\right).$$ It is clear that
$\det(BC_{t,s}-C_{t,s}B)\not=0$ for some $t,s$ whenever
$\alpha\not\in{\mathbb R}$ or $k\alpha\not=\frac{\bar{\alpha}}{k}$
or $k\alpha\not=\xi_3$ or $\xi_3\not=\frac{\bar{\alpha}}{k}$, and in
this case we have $h(B)=h(C_{t,s})=h(A)$. If $\alpha$ is real and
$k\alpha=\frac{\bar{\alpha}}{k}=\xi_3$, then, up to a real scalar
multiple, $B$ has the form
$$B=\left(\begin{array}{ccc} 1 &1 &0 \\ 1 &1 & 0\\
0&0&1\end{array}\right).$$ Let
$$C=\left(\begin{array}{ccc} 1 &1 &1+i \\ 1 &2 & 1-i\\
1-i&1+i&0\end{array}\right).$$ Then ${\rm Im}(1\cdot
(1-i)\overline{(1+i)})=-2i\not=0$ and hence $h(C)=h(A)$. Since
$\det(BC-CB)=-4i\not=0$, we also have $h(B)=h(C)$. So, again we get
$h(B)=h(A)$, as desired.

{\bf Subcase 2.} One of  $\alpha,\beta,\gamma$ is 0.

Without loss of generality, say $\beta=0$. Then, as rank $B=2$,
$\det B=\xi_1\xi_2\xi_3-|\gamma|^2\xi_2-|\alpha|^2\xi_3=0$. Thus
there are scalars $c,d$ with $d\not=0$ such that
$c\xi_1=\bar{\gamma}-d\alpha$, $\xi_2=-\frac{c}{d}\alpha$ and
$\xi_3=c\gamma$.

Clearly, $\xi_2=0\Leftrightarrow\xi_3=0\Leftrightarrow c=0$, and in
this case we have
$$B=\left(\begin{array}{ccc} \xi_{1} &\alpha &\gamma \\ \bar{\alpha} &0& 0\\
\bar{\gamma}&0&0\end{array}\right).$$ Let
\begin{equation} C_{t,s,p}=\left(\begin{array}{ccc} 0 &t & ip \\ t &0 & s\\
-ip&s&0\end{array}\right)\end{equation} for nonzero real numbers
$t,s,p$; then $h(A)=h(C_{t,s,p})$. Now
$$BC_{t,s,p}-C_{t,s,p}B=\left(\begin{array}{ccc} t(\alpha-\bar{\alpha})-ip(\gamma+\bar{\gamma}) &\xi_1t+\gamma s &i\xi_1p+\alpha s \\
-\xi_1t-\bar{\gamma}s  &-t(\alpha-\bar{\alpha})& i\bar{\alpha}p-\gamma t\\
i\xi_1p-\bar{\alpha}s&t\bar{\gamma}+i\alpha
p&ip(\gamma+\bar{\gamma})\end{array}\right).$$ If $\xi_1\not=0$ (in
this case the coefficients of $sp^2$ and $t^2s$ of
$\det(BC_{t,s,p}-C_{t,s,p}B)$ are nonzero), or if $\xi_1=0$ but one
of $\alpha-\bar{\alpha}$ and $\gamma+\bar{\gamma}$ is nonzero (in
this case the coefficient of $t^3$ or $p^3$ is nonzero), it is sure
that $BC_{t,s,p}-C_{t,s,p}B$ is of rank three for some $t,s,p$ and
hence $h(B)=h(C_{t,s,p})=h(A)$. If
$$B=\left(\begin{array}{ccc} 0 &\alpha &i\delta \\ {\alpha} &0& 0\\
-i\delta &0&0\end{array}\right) $$ for some nonzero real numbers
$\alpha,\delta$, let
\begin{equation} D_{t,s,p}=\left(\begin{array}{ccc} 0 &it & p \\ -it &0 & s\\
p&s&0\end{array}\right)\end{equation} for nonzero real numbers
$t,s,p$. Then
$$BD_{t,s,p}-D_{t,s,p}B=\left(\begin{array}{ccc}2i(\delta p-\alpha t) &i\delta s & \alpha s \\ i\delta s &2i\alpha t & \alpha p-\delta t\\
-\alpha s&\delta t-\alpha p&-2i\delta p\end{array}\right),$$ which
is of rank three for some suitable choice of $t,s,p$ as the
coefficients of $t^3$ and $p^3$ of $\det(BD_{t,s,p}-D_{t,s,p}B)$ are
nonzero. Therefore, we have $h(B)=h(D_{t,s,p})=h(A)$.

Assume that $c\not=0$; then
$\xi_1=\frac{1}{c}(\bar{\gamma}-d\bar{\alpha}),
\xi_2=-\frac{c}{d}\alpha, \xi_3=c\gamma$ are real,
$$ B=\left(\begin{array}{ccc} \frac{1}{c}(\bar{\gamma}-d\bar{\alpha}) &\alpha &\gamma   \\ \bar{\alpha} &-\frac{c}{d}\alpha & 0\\
\bar{\gamma}&0&c\gamma\end{array}\right)$$ and, for $C_{t,s,p}$ in
Eq.(3.5), we have
$$\begin{array}{rl}&BC_{t,s,p}-C_{t,s,p}B\\=&
\left(\begin{array}{ccc} t(\alpha-\bar{\alpha})-ip(\gamma+\bar{\gamma}) &(\frac{\bar{\gamma}}{c}-\frac{d\bar{\alpha}}{c}+\frac{c\alpha}{d})t+\gamma s
&i(\frac{\bar{\gamma}}{c}-\frac{d\bar{\alpha}}{c}-c\gamma)p+\alpha s \\
-(\frac{\bar{\gamma}}{c}-\frac{d\bar{\alpha}}{c}+\frac{c\alpha}{d})t-\bar{\gamma}s  &-t(\alpha-\bar{\alpha})& i\bar{\alpha}p-\gamma t-c(\frac{\alpha}{d}+\gamma)s\\
i(\frac{\bar{\gamma}}{c}-\frac{d\bar{\alpha}}{c}-c\gamma)p-\bar{\alpha}s&i\alpha
p+t\bar{\gamma}+c(\frac{\alpha}{d}+\gamma)s&ip(\gamma+\bar{\gamma})\end{array}\right)\end{array}.$$
Note that the coefficients of $t^3,s^3$ and $p^3$ in
$\det(BC_{t,s,p}-C_{t,s,p}B)$ are respectively
$|\gamma|^2(\alpha-\bar{\alpha}),
c(\frac{\alpha}{d}+\gamma)(\bar{\alpha}\gamma-\alpha\bar{\gamma})$
and $-i|\alpha|^2(\gamma+\bar{\gamma})$.

It is clear that, if $\alpha$ or $i\gamma$ are not real; or in the
case that both $\alpha$ and $i\gamma$ are   real,  but $
\xi_2\not=\xi_3$, then $BC_{t,s,p}-C_{t,s,p}B$ is rank-3 for
suitable choice of real numbers $t,s,p$ and hence
$h(B)=h(C_{t,s,p})=h(A)$.

If $\alpha,i\gamma$ are real and $\xi_2=\xi_3$ but
$\xi_1\not=\xi_2$, then
$$\begin{array}{rl}&BC_{t,s,p}-C_{t,s,p}B\\=&
\left(\begin{array}{ccc} 0
&(\frac{\bar{\gamma}}{c}-\frac{d\bar{\alpha}}{c}+\frac{c\alpha}{d})t+\gamma
s
&i(\frac{\bar{\gamma}}{c}-\frac{d\bar{\alpha}}{c}-c\gamma)p+\alpha s \\
-(\frac{\bar{\gamma}}{c}-\frac{d\bar{\alpha}}{c}+\frac{c\alpha}{d})t-\bar{\gamma}s  &0& i\bar{\alpha}p-\gamma t \\
i(\frac{\bar{\gamma}}{c}-\frac{d\bar{\alpha}}{c}-c\gamma)p-\bar{\alpha}s&i\alpha
p+t\bar{\gamma}&0\end{array}\right)\end{array}.$$ As the coefficient
of $t^2p$ in $\det(BC_{t,s,p}-C_{t,s,p}B)$ is
$-i(\xi_1-\xi_2)^2\bar{\gamma}\not=0$, we still have $h(B)=h(A)$.

If $\alpha,i\gamma$ are real and $\xi_1=\xi_2=\xi_3$, then $B$ has
the form
$$ B=\left(\begin{array}{ccc} \pm\sqrt{\alpha^2+\delta^2} &\alpha &i\delta   \\ {\alpha} &\pm\sqrt{\alpha^2+\delta^2}& 0\\
-i\delta&0&\pm\sqrt{\alpha^2+\delta^2}\end{array}\right)$$ with
nonzero $\alpha, \delta\in{\mathbb R}$. Then, for $D_{t,s,p}$ in
Eq.(3.6), consider
$$BD_{t,s,p}-D_{t,s,p}B=\left(\begin{array}{ccc}2i(\delta p-\alpha t) &i\delta s & \alpha s \\ i\delta s &2i\alpha t & \alpha p-\delta t\\
-\alpha s&\delta t-\alpha p&-2i\delta p\end{array}\right)$$ for
nonzero real numbers $t,s,p$. As the coefficient of $t^3$ in
$\det(BD_{t,s,p}-D_{t,s,p}B)$ is $-2i\alpha\delta^2\not=0$, and
hence one gets $h(B)=h(A)$ again.

{\bf Subcase 3.} All $\alpha,\beta,\gamma$ are nonzero.

Since $\det(B)=0$, there are scalars $c,d$ such that $(\bar{\gamma},
\bar{\beta}, \xi_3)=(c\xi_1+d\bar{\alpha},
c\alpha+d\xi_2,c\gamma+d\beta)$. It follows that
\begin{equation}\left\{\begin{array}{l} \gamma=\bar{c}\xi_1+\bar{d}\alpha ,\\
\beta=\bar{c}\bar{\alpha}+\bar{d}\xi_2, \\
\xi_3=|c|^2\xi_1+c\bar{d}\alpha+\bar{c}d\bar{\alpha}+|d|^2\xi_2.\end{array}\right.
\end{equation}
As $\alpha\beta\bar{\gamma}\in{\mathbb R}$, we get
$$(c\bar{d}\alpha-\bar{c}d\bar{\alpha})(|\alpha|^2-\xi_1\xi_2)=0.$$
However, $|\alpha|^2-\xi_1\xi_2=0$ implies that
$\xi_1=\frac{\bar{\alpha}}{k}$, $\xi_2=k\alpha$ for some scalar $k$,
which entails that $\beta=k\gamma$ and hence $B$ is of rank-1, a
contradiction. So $|\alpha|^2-\xi_1\xi_2\not=0$ and then we must
have $c\bar{d}\alpha-\bar{c}d\bar{\alpha}=0$. Discussing similarly,
we get
\begin{equation} \left\{ \begin{array}{l}
|\alpha|^2-\xi_1\xi_2\not=0,\\ |\beta|^2-\xi_2\xi_3\not=0,\\
 |\gamma|^2-\xi_1\xi_3\not=0. \end{array} \right. \end{equation}

 Let $C_{t,s,p}$ be as in Eq.(3.5). As
$$\begin{array}{rl} &BC_{t,s,p}-C_{t,s,p}B\\ =& \left(\begin{array}{ccc}
(\alpha-\bar{\alpha})t-i(\gamma+\bar{\gamma})p &
(\xi_1-\xi_2)t+\gamma s-i\bar{\beta}p & i(\xi_1-\xi_3)p+\alpha
s-\beta t \\ -(\xi_1-\xi_2)t-\bar{\gamma}s-i\beta p &
-(\alpha-\bar{\alpha})t+(\beta-\bar{\beta})s &
(\xi_2-\xi_3)s+i\bar{\alpha}p-\gamma t \\
i(\xi_1-\xi_3)p+\bar{\beta}t-\bar{\alpha}s &
-(\xi_2-\xi_3)s+\bar{\gamma}t+i\alpha p &
i(\gamma+\bar{\gamma})p-(\beta-\bar{\beta})s
\end{array}\right),\end{array}$$
we see that the coefficients of $t^3, s^3, p^3$ in
$\det(BC_{t,s,p}-C_{t,s,p}B)$ are respectively
\begin{equation}\left\{ \begin{array}{l}
c_t=(\xi_1-\xi_2)(\beta\bar{\gamma}-\bar{\beta}\gamma)+(\alpha-\bar{\alpha})(|\gamma|^2-|\beta|^2),\\
c_s=(\xi_2-\xi_3)(\alpha\bar{\gamma}-\bar{\alpha}\gamma)+(\beta-\bar{\beta})(|\alpha|^2-|\gamma|^2),\\
c_p=i(\xi_1-\xi_3)(\bar{\alpha}\bar{\beta}+\alpha\beta)+i(\bar{\gamma}+\gamma)(|\beta|^2-|\alpha|^2).
\end{array}\right. \end{equation}
 If one of  $c_t, c_s,c_p$ is nonzero, then
$\det(BC_{t,s,p}-C_{t,s,p}B)\not=0$ for some choice of $t,s,p$,
which implies that $h(B)=h(C_{t,s,p})=h(A)$. Assume
$$c_t=c_s=c_p=0.$$ Considering the coefficients $d_t,d_s$ and $d_p$
of $t^3,s^3$ and $p^3$ in $\det(BD_{t,s,p}-D_{t,s,p}B)$  with
$D_{t,s,p}$ as in Eq.(3.6) one gets
\begin{equation}
\left\{ \begin{array}{l}
d_t=i(\xi_1-\xi_2)(\beta \bar{\gamma}+\bar{\beta}{\gamma})+i(\alpha+\bar{\alpha})(|\gamma|^2-|\beta|^2),\\
d_s=(\xi_2-\xi_3)(\alpha\bar{\gamma}-\bar{\alpha}\gamma)+(\beta-\bar{\beta})(|\alpha|^2-|\gamma|^2)=0,\\
d_p=(\xi_1-\xi_3)(\bar{\alpha}\bar{\beta}-
\alpha{\beta})+(\bar{\gamma}-{\gamma})(|\beta|^2-|\alpha|^2).
\end{array}\right.
\end{equation}
If one of $d_t,d_p$ is nonzero, then $h(B)=h(A)$. Assume that
$$d_t=d_s=d_p=0.$$
Let
\begin{equation} E_{t,s,p}=\left(\begin{array}{ccc} 0 &t & p \\ t &0 & is\\
p&-is&0\end{array}\right)\end{equation} for nonzero real numbers
$t,s,p$. The coefficients $e_t,e_s$ and $e_p$ of $t^3,s^3$ and $p^3$
in $\det(BE_{t,s,p}-E_{t,s,p}B)$ are
\begin{equation}
\left\{ \begin{array}{l}
e_t= (\xi_1-\xi_2)(\beta \bar{\gamma}-\bar{\beta}{\gamma})+ (\alpha-\bar{\alpha})(|\gamma|^2-|\beta|^2)=0,\\
e_s=i(\xi_2-\xi_3)(\alpha\bar{\gamma}+\bar{\alpha}\gamma)+i(\beta+\bar{\beta})(|\alpha|^2-|\gamma|^2),\\
e_p=(\xi_1-\xi_3)(\bar{\alpha}\bar{\beta}-
\alpha{\beta})+(\bar{\gamma}-{\gamma})(|\beta|^2-|\alpha|^2)=0.
\end{array}\right.
\end{equation}
If $e_s\not=0$, then we get $h(B)=h(A)$. Assume
$$e_s=0.$$
Then, by Eqs.(3.9)-(3.10), and Eq.(3.12), it is easily checked that
\begin{equation}
\left\{ \begin{array}{l}
(\xi_1-\xi_2)\beta \bar{\gamma}+\alpha(|\gamma|^2-|\beta|^2)=0,\\
(\xi_2-\xi_3)\alpha\bar{\gamma}+\beta(|\alpha|^2-|\gamma|^2)=0,\\
(\xi_1-\xi_3) \alpha{\beta}+{\gamma}(|\beta|^2-|\alpha|^2)=0.
\end{array}\right.
\end{equation}

As $\alpha\beta\bar{\gamma}$ is real, we see from Eq.(3.13) that
both $\alpha^2, \beta^2$ are real and hence $\alpha\in{\mathbb R}$
or $\alpha\in i{\mathbb R}$ ($\beta\in{\mathbb R}$ or $\beta\in
i{\mathbb R}$). It follows that there are four cases may occur, that
is,
\begin{equation}\left\{\begin{array}{ll}
1^\circ & \alpha,\beta,\gamma\in{\mathbb R}.\\
2^\circ & \alpha,\beta\in i{\mathbb R},\gamma\in{\mathbb R}.\\
3^\circ
& \beta,\gamma\in i{\mathbb R}, \alpha\in {\mathbb R}.\\
4^\circ & \alpha,\gamma\in i{\mathbb R},\beta\in{\mathbb R}.
\end{array}\right.\end{equation}

If $\xi_1=\xi_2=\xi_3=\xi$, then $|\alpha|=|\beta|=|\gamma|$ by
Eq.(3.13). On the other hand, by Eq.(3.7),
$\xi(1-|c|^2-|d|^2)=2c\bar{d}\alpha$. Thus $\xi=0$ or
$1-|c|^2-|d|^2=0$ implies that $c=0$ or $d=0$. Without loss of
generality, say $c=0$; then $d\not=0$ and $\xi=|d|^2\xi\not=0$,
which gives $d=e^{i\theta}$ and $|\xi|=|d\beta|=|\alpha|$,
contradicting the fact that
$|\alpha|^2\not=\xi_1\xi_2=\xi^2=|\xi|^2$ (see Eq.(3.8)).

  So we have $\xi(1-|c|^2-|d|^2)=2c\bar{d}\alpha\not=0$,
$$\xi=\frac{\gamma-\bar{d}\alpha}{c}=\frac{\beta-\bar{c}\bar{\alpha}}{\bar{d}}=\frac{2c\bar{d}\alpha}{1-|c|^2-|d|^2}=\frac{2\bar{c}d\bar{
\alpha}}{1-|c|^2-|d|^2}$$ and
$$B=\left(\begin{array}{ccc} \frac{2c\bar{d}\alpha}{1-|c|^2-|d|^2} &\alpha &(2|c|^2+1)\bar{d}\alpha  \\ \bar{\alpha} &\frac{2c\bar{d}\alpha}{1-|c|^2-|d|^2}& (2|d|^2+1)\bar{c}\bar{\alpha}\\
(2|c|^2+1)d\bar{\alpha}&(2|d|^2+1){c}{\alpha}&\frac{2c\bar{d}\alpha}{1-|c|^2-|d|^2}\end{array}\right).$$
It follows that $(2|c|^2+1)|d|=(2|d|^2+1)|c|=1$ as
$|\alpha|=|\beta|=|\gamma|$. Thus
$2|c|+\frac{1}{|c|}=2|d|+\frac{1}{|d|}=\frac{1}{|cd|}$. Note that
$|c|=\frac{1}{2|d|^2+1}$ and $|d|=\frac{1}{2|c|^2+1}$. So one gets
$(|c|+|d|)(|c|-|d|)=|c|-|d|$, which gives further that either
$|c|=|d|$ or $|c|+|d|=1$. If $|c|\not=|d|$, we must have $|c|+|d|=1$
and hence $0<1-|c|=|d|=\frac{1}{2|c|^2+1}$. Then we obtain
$|c|(2|c|^2-2|c|+1)=0$. As we always have $2|c|^2-2|c|+1>0$, one
sees that $c=0$, a contradiction.  Therefore, we have $|c|=|d|=k$.
Since $(2k^2+1)k=1$, we see that $k\approx 0.5898$. Write
$\alpha=|\alpha|e^{i\theta_1}$, $c=ke^{i\theta_2}$ and
$d=ke^{i\theta_3}$. Now $c\bar{d}\alpha$ is real implies that
$\theta_1+\theta_2-\theta_3$ is 0 or $\pi$. Replacing $B$ by $-B$ if
necessary we may assume that $\theta_1+\theta_2-\theta_3=0$ and thus
$d=ke^{i(\theta_1+\theta_2)}$. Without loss of generality, let
$|\alpha|=1$. Notice that $(2k^2+1)k=1$. Then $B$ becomes to
$$B=\left(\begin{array}{ccc} \frac{2k^2}{1-2k^2} &e^{i\theta_1} & e^{-i\theta_2}\\ e^{-i\theta_1} &\frac{2k^2}{1-2k^2}& e^{-i(\theta_1+\theta_2)}\\
 e^{i\theta_2}& e^{i(\theta_1+\theta_2)}&\frac{2k^2}{1-2k^2}\end{array}\right)$$
 with $\frac{2k^2}{1-2k^2}\approx 2.2868$. But then
 $0=\det(B)=(\frac{2k^2}{1-2k^2})^3-3(\frac{2k^2}{1-2k^2})+2\approx
 7.0983>0$, a contradiction. Therefore $\xi_1,\xi_2, \xi_3$ are not all the
 same. Keep this in mind below, we can show that $h(B)=h(A)$ holds.

 For example, consider the Case 2$^\circ$,  that is, $\alpha\in{\mathbb R},\beta,\gamma \in i{\mathbb R}$.

In this case, for $t,s,p\in{\mathbb C}$, let
$$F_{t,s,p}=\left(\begin{array}{ccc} 0 &t & p \\ \bar{t} & 0 & s
\\ \bar{p} & \bar{s} & 0 \end{array}\right).$$
 Then
$$\small\begin{array}{rl} &BF_{t,s,p}-F_{t,s,p}B \\
= &\left(\begin{array}{ccc} \alpha(\bar{t}-t)+i\gamma(p+\bar{p}) &  (\xi_1-\xi_2)t+i\gamma\bar{s}+i\beta p  &  (\xi_1-\xi_3)p-\alpha s-i\beta t \\
(\xi_2-\xi_1)\bar{t}+i\beta\bar{p}+i\gamma  s &
-\alpha(\bar{t}-t)+i\beta (s+\bar{s})& (\xi_2-\xi_3)s+\alpha
p-i\gamma\bar{t}
\\ (\xi_3-\xi_1)\bar{p}-i\beta\bar{t}-\alpha\bar{s} & (\xi_3-\xi_2)\bar{s}-i\gamma t-\alpha\bar{p} & -i\gamma (p+\bar{p})-i\beta(s+\bar{s})
 \end{array}\right).\end{array}$$
Consider the term of $\det([B,F_{t,s,p}])$  that  contains only $t$,
which is
$$
((\xi_1-\xi_2)\beta\gamma
+\alpha(\beta^2-\gamma^2))(t^2\bar{t}-\bar{t}^2t)=2(\xi_1-\xi_2)\beta\gamma(t^2\bar{t}-\bar{t}^2t)
$$
as
$(\xi_1-\xi_2)\beta\gamma+\alpha(\gamma^2-\beta^2)=(\xi_1-\xi_2)(i\beta)\bar{i\gamma}+\alpha(|i\gamma|^2-|i\beta|^2)=0$
by Eq.(3.13). If $\xi_1\not=\xi_2$, then
$(\xi_1-\xi_2)\beta\gamma\not=0$ and it is clear that we can choose
$t,s,p$ with $ts\bar{p}\not\in{\mathbb R}$ so that
$\det([B,F_{t,s,p}])\not=0$. Thus we get $h(B)=h(F_{t,s,p})=h(A)$.
If $\xi_1=\xi_2$, then we must have $\xi_2\not=\xi_3$. Now consider
the term of $\det([B,F_{t,s,p}])$ that only contain $s$, which is
$$
i(\xi_2-\xi_3)\alpha\gamma
(s^2\bar{s}-\bar{s}^2s)+i\beta(\alpha^2-\gamma^2)(s^2\bar{s}+\bar{s}^2s)=2i(\xi_2-\xi_3)\alpha\gamma
s^2\bar{s}
$$
since
$(\xi_2-\xi_3)\alpha\bar{i\gamma}+(i\beta)(|\alpha|^2-|i\gamma|^2)=0$
by Eq.(3.13). Clearly  $(\xi_2-\xi_3)\alpha\gamma\not=0$  implies
that there are $t,s,p$ with $ts\bar{p}\not\in{\mathbb R}$ so that
$\det([B,F_{t,s,p}])\not=0$. It follows that
$h(B)=h(F_{t,s,p})=h(A)$.

The cases 1$^\circ$, 3$^\circ$ and 4$^\circ$ are dealt with
similarly. This completes the proof of the Claim 4.

{\bf Claim 5.} For any $A, B\in{\mathcal B}_s(H)\setminus{\mathcal
D}$, we have $h(A)=h(B)$.

By Lemma 3.2, there exist $E,F\in {\mathcal
B}_s(H)\setminus{\mathcal D}$ of rank not greater than 2 such that
$W(AE-EF)\not=-W(AE-EF)$ and $W(BF-FB)\not=-W(BF-FB)$. Thus we get
$h(A)=h(E)$ and $h(B)=h(F)$. However, by Claims 2-4, we always have
$h(E)=h(F)$. Hence $h(A)=h(B)$.

Finally, let ${\mathcal S}=\{S\in{\mathcal D}: h(S)\not=h(A) \ {\rm
for}\ A\not\in{\mathcal D}\}$. Then it is clear that the theorem
holds. \hfill$\Box$

\section{The case when $\dim H=2$}

In this last section we consider the problem for the case when $\dim
H=2$. As we will see the situation for the two dimensional case is
much different from that for the case of dimension $\geq 3$.

As $\dim H=2$, we can identify ${\mathcal B}_s(H)$ as ${\bf
H}_2={\bf H}_2({\mathbb C})$, the set of all $2\times 2$ Hermitian
matrices over $\mathbb C$.

The following is our result and the  surjectivity assumption on
$\Phi$ is not needed.

{\bf Theorem 4.1.} {\it Let $\Phi: {\bf H}_2({\mathbb C})\to{\bf
H}_2({\mathbb C})$ be a map. The following statements are
equivalent.}

(1) {\it $\sigma([\Phi(A),\Phi(B)])=\sigma([A,B])$ for any
$A,B\in{\bf H}_2({\mathbb C})$.}

(2) {\it $W([\Phi(A),\Phi(B)])=W([A,B])$ for any $A,B\in{\bf
H}_2({\mathbb C})$.}

(3) {\it $w([\Phi(A),\Phi(B)])=w([A,B])$ for any $A,B\in{\bf
H}_2({\mathbb C})$.}

(4) {\it There exist a unitary matrix $U\in M_2({\mathbb C})$, a
sign function $h:{\bf H}_2\to \{-1,1\}$ and a functional $f:{\bf
H}_2({\mathbb C})\to {\mathbb R}$ such that   one of the following
holds:}

\hspace{2mm}  (1$^\circ$) {\it  $\Phi(A)=h(A)UAU^*+f(A)I$ for all
$A\in{\bf H}_2$;}

\hspace{2mm}  (2$^\circ$) {\it  $\Phi(A)=h(A)UA^tU^*+f(A)I$ for all
$A\in{\bf H}_2$;}

\hspace{2mm} (3$^\circ$) {\it $\Phi(A)=h(A)U\Psi(A)U^*+f(A)I$ for
all $A\in{\bf H}_2$;}

\hspace{2mm} (4$^\circ$) {\it $\Phi(A)=h(A)U\Psi(A)^tU^*+f(A)I$ for
all $A\in{\bf H}_2$.}\\
Where, with $A= \left(\begin{array}{cc} a & c+id \\ c-id &b
\end{array}\right)$, $\Psi(A)=\left(\begin{array}{cc} a & -c+id \\ -c-id &b
\end{array}\right)$.
 \vspace{3mm}

 {\bf Proof.} It is clear that
 (4)$\Rightarrow$(1)$\Leftrightarrow$(2)$\Leftrightarrow$(3).

 (3)$\Rightarrow$(4). Assume $\Phi:{\bf H}_2\to{\bf H}_2$ preserves
 the numerical radius of Lie product.

 We may modify the
functional $f(A)$ in the map $\Phi$ so that $\Phi(A)$ has trace 0
for all $A \in {\bf H}_2({\mathbb C})$. Then we can focus on the set
${\bf H}_2^0$ of trace zero matrices in ${\bf H}_2({\mathbb C})$.

Now, suppose (1) holds.

Consider the Hermitian matrices
\begin{equation} X = \frac{1}{\sqrt
2} \left(\begin{array}{cc} 0 & 1 \cr 1 & 0 \cr
\end{array}\right),  \qquad
Y = \frac{1}{\sqrt 2}\left(\begin{array}{cc} 0 & -i \cr i & 0 \cr
\end{array}\right),  \qquad
Z = \frac{1}{\sqrt 2}\left(\begin{array}{cc} 1 & 0 \cr 0 & -1 \cr
\end{array}\right).
\end{equation}

Then the following holds:

(1)  $\{X, Y, Z\}$ is an orthonormal  basis for $M_2^0$ using the
inner product $\langle A, B \rangle = {\rm tr}( AB^*)$, where
$M_2^0$ is the set of trace zero $2\times 2$ matrices.

(2)   $A = a_1 X + a_2 Y + a_3 Z\in{\bf H}_2^0$ if and only if
$(a_1, a_2, a_3)^t\in{\mathbb R}^3$.

(3) $XY = \frac{i}{\sqrt 2}Z = -YX, \quad YZ = \frac{i}{\sqrt 2}X =
-ZY, \quad ZX = \frac{i}{\sqrt 2}Y = -XZ$.

(4) $W([X,Y]) = W([Y,Z]) = W([Z,X]) = i[-1,1]$.

(5)  If $A = a_1 X + a_2 Y + a_3 Z$ and $B = b_1 X + b_2 Y + b_3 Z$
in $M_2^0$, then
$$[A,B] = \sqrt{2}i(c_1 X + c_2 Y + c_3 Z),$$
where
$$c_1 = a_2 b_3 - a_3b_2, \quad c_2 = -(a_1 b_3 - a_3 b_1), \quad c_3 = a_1 b_2 - a_2 b_1.$$
In other words, $(c_1, c_2, c_3)^t = (a_1, a_2, a_3)^t \times (b_1,
b_2, b_3)^t$, the cross product in ${\mathbb C}^3$.

(6) Every unitary similarity map $a_1X + a_2Y + a_3Z = A \mapsto
UAU^* = b_1X + b_2 Y + b_3Z$ on  $M_2^0$ corresponds to a real
special orthogonal transformation $T \in M_3({\mathbb C})$ such that
$T(a_1, a_2, a_3)^t = (b_1, b_2, b_3)^t$.

{\bf Claim 1.}  There exist a unitary $U\in M_2({\mathbb C})$  such
that
$$\Phi(A) = \varepsilon_A UAU^*$$ for all $A \in \{X, Y, Z\},$
where $\varepsilon_A\in\{-1,1\}$.

 Assume that the image of $X, Y, Z$ are respectively
$$X_1 = a_{11} X + a_{21} Y + a_{31} Z,
\  Y_1 = a_{12} X + a_{22} Y + a_{32} Z, \  Z_1 = a_{13} X + a_{23}
Y + a_{33} Z.$$   Then $a_{pq}$s are real numbers. Let $T = (a_{pq})
\in M_3({\mathbb R})$. We will show that   $T$ is a real orthogonal
matrix. Thus $\Phi$ has the form in Claim 6.

Note that the hypothesis and conclusion will not be affected by
changing $T$ to $PTQ$ for any real orthogonal matrices $P, Q \in
M_3({\mathbb C})$. It just corresponds to changing $\Phi$ to a map
of the form
$$A \mapsto \varepsilon_P U_P\Phi(\varepsilon_Q U_Q A U_Q^*)U_P^*$$
for some unitary $U_P, U_Q \in M_2({\mathbb C})$ and $\varepsilon_P,
\varepsilon_Q \in \{1,-1\}$ depending on $P$ and $Q$.

 By the
singular value decomposition of real matrices, let $P, Q$ be real
orthogonal such that $PTQ = {\rm diag}(s_1, s_2, s_3)$ with $s_1 \ge
s_2 \ge s_3 \ge 0$. Now, replace $T$ by $PTQ$ so that   $T = {\rm
diag}(s_1, s_2, s_3)$. Thus there exists a real orthogonal matrix
$U\in M_2({\mathbb C})$ such that
$$ \Phi(X)= s_1UXU^*,\  \Phi(Y) =s_2UYU^*,\  \Phi(Z) =s_3UZU^*.$$
It follows that
$$ 1=w(XY-YX)=w( \Phi(X)\Phi(Y)- \Phi(Y)\Phi(X) )=|s_1s_2|w(XY-YX)=|s_1s_2|.$$
Similarly, one gets $|s_1s_3|=|s_2s_3|=1$ and  hence $s_1,s_2,s_3
\in\{-1,1\}$. Thus the Claim is true.

Without loss of generality in the sequel we assume $U=I_2$. Note
that, for any sign function $h: {\bf H}_2\to\{-1,1\}$, the map
$\Psi$ defined by $\Psi (A)=h(A)\Phi(A)$ still preserves the
numerical radius of Lie product. So, multiplied by a suitable sign
function if necessary, we may assume that
$$\Phi(C)=C$$ for every
$C\in\{X,Y,Z\}$.

{\bf Claim 2.} There are sign functions
$\varepsilon_1,\varepsilon_2,\varepsilon_3:{\bf H}_2\to\{-1,1\}$ and
a functional $f:{\bf H}_2\to{\mathbb R}$ such that, for
any $A \in {\bf H}_2$ with $A=\left(\begin{array}{cc} a & c+id \\
c-id &b
\end{array}\right)$, we have $$\Phi(A)=\left(\begin{array}{cc} \varepsilon_1(A)a & \varepsilon_2(A)c+i\varepsilon_3(A)d \\ \varepsilon_2(A)c+i\varepsilon_3(A)d &
\varepsilon_1(A)b
\end{array}\right)+f(A)I_2.$$

Write $A=\left(\begin{array}{cc} a & c+id \\ c-id &b
\end{array}\right)$ and $\Phi(A)=\left(\begin{array}{cc} x & w+iv \\ w-iv
&y
\end{array}\right)$, where $a,b,c,d,x,y,w,v$ are real numbers. Note that, for any $E,F\in {\bf H}_2$,
$w(EF-FE)=\delta$ if and only if the spectrum
$\sigma(EF-FE)=i[-\delta,\delta]$. Thus $w(AC-CA)=w(BC-CB)$ if and
only if $\sigma(AC-CA)=\sigma(BC-CB)$. As
$$\begin{array}{ll} \sqrt{2}(AX-XA)= \left(\begin{array}{cc} i2d & a-b \\ b-a &-i2d
\end{array}\right), & \sqrt{2}(\Phi(A)X-X\Phi(A))= \left(\begin{array}{cc} 2iv & x-y \\
y-x &-2iv
\end{array}\right);\\ \sqrt{2}(AY-YA)=i\left(\begin{array}{cc} 2c & b-a \\ b-a
&-2c
\end{array}\right), & \sqrt{2}(\Phi(A)Y-Y\Phi(A))=i\left(\begin{array}{cc} 2w & y-x \\
y-x&-2w
\end{array}\right); \\ \sqrt{2}(AZ-ZA)=2\left(\begin{array}{cc} 0 & -c-id
\\c-id
&0
\end{array}\right),& \sqrt{2}(\Phi(A)Z-Z\Phi(A))=2\left(\begin{array}{cc} 0 &
-w-iv
\\w-iv
&0
\end{array}\right),\end{array}$$
we must have \begin{equation}
\left\{\begin{array}{l} 4v^2+(x-y)^2=4d^2+(a-b)^2, \\
4w^2-(x-y)^2=4c^2-(a-b)^2, \\w^2+v^2=c^2+d^2.\end{array}\right.
\end{equation}

It follows that, the map $\Phi$ sends $\left(\begin{array}{cc} a & 0
\\ 0 &b
\end{array}\right)=\left(\begin{array}{cc} \frac{a-b}{2} & 0
\\ 0 &-\frac{a-b}{2}
\end{array}\right)+\frac{a+b}{2}I_2$ to $\left(\begin{array}{cc} \varepsilon_1 \frac{a-b}{2} & 0\\ 0
& -\varepsilon_1 \frac{a-b}{2}
\end{array}\right)+\lambda' I_2=\left(\begin{array}{cc} \varepsilon_1 a & 0\\ 0
& \varepsilon_1 b
\end{array}\right)+\lambda I_2$,  and sends $\left(\begin{array}{cc} 0 & c+id \\ c-id
&0
\end{array}\right)$ to $\left(\begin{array}{cc} 0 & \varepsilon_2c+i\varepsilon_3d \\ \varepsilon_2c-i\varepsilon_3d
&0
\end{array}\right)+\lambda_2I_2$ for some scalars
$\varepsilon_1,\varepsilon_2,\varepsilon_3\in\{-1,1\}$.

To sum up,
$$\Phi(\left(\begin{array}{cc} a & 0
\\ 0 &b
\end{array}\right)+{\mathbb R}I_2)\subseteq\varepsilon_1\left(\begin{array}{cc} a & 0
\\ 0 &b
\end{array}\right)+{\mathbb R}I_2,$$ and
$$\Phi(\left(\begin{array}{cc} 0 & c+id
\\ c-id &0
\end{array}\right)+{\mathbb R}I_2)\subseteq \left(\begin{array}{cc} 0 & \varepsilon_2c+i\varepsilon_3d \\ \varepsilon_2c-i\varepsilon_3d
&0
\end{array}\right)+{\mathbb R}I_2,$$ where
$\varepsilon_1,\varepsilon_2,\varepsilon_3\in\{-1,1\}$ depending on
$a,c,d$.

To consider the general $A=\left(\begin{array}{cc} a & c+id \\
c-id &b
\end{array}\right)$, for any unit vector $x\in{\mathbb C}^2$,
take unit vector $y\perp x$. Then, with respect to the orthonormal
base $\{x,y\}$, one can take $$X'=\frac{1}{\sqrt{2}}(x\otimes
y+y\otimes x),\ Y'=\frac{1}{\sqrt{2}}i(-x\otimes y+y\otimes x),\
Z'=\frac{1}{\sqrt{2}}(x\otimes x-y\otimes y).$$ Repeat the argument
as in Claim 1 and the above one achieves that, there exists a
unitary matrix $U_x$ such that
\begin{equation}
\Phi(ax\otimes x+by\otimes y+{\mathbb R}I_2)\subseteq
\varepsilon_1(x,a,b)(aU_xx\otimes U_x x  +bU_xy\otimes U_xy )
+{\mathbb R}I_2
\end{equation}
  for any $a,b\in{\mathbb R}$ and
\begin{equation}
\begin{array}{rl} &\Phi((c+id)x\otimes y+(c-id)y\otimes x+{\mathbb R}I_2) \\
\subseteq &
(\varepsilon_2(x,c,d)c+i\varepsilon_3(x,c,d)d)u_xx\otimes U_xy
\\ &+(\varepsilon_2(x,c,d)c-i\varepsilon_3(x,c,d)d)U_xy\otimes
U_xx+{\mathbb R}I_2, \end{array}
\end{equation}
where
$\varepsilon_1(x,a,b),\varepsilon_2(x,c,d),\varepsilon_3(x,c,d))\in\{-1,1\}$.
Particularly, by Eq.(4.3), without loss of generality we may assume
that \begin{equation}\sigma(\Phi(A))=\sigma(A) \end{equation} for
all $A\in {\bf H}_2$. It follows that, if $b=-a$, that is, if
$A\in{\bf H}_2^0$, then we have
\begin{equation}
x^2+w^2+v^2=a^2+c^2+d^2,
\end{equation}
which, together with Eq.(4.2), gives
$$ x^2=a^2,\ w^2=c^2,\ v^2=d^2.$$
Therefore, we still have
$$x=\varepsilon_1 a,\ w=\varepsilon_2c, \ v=\varepsilon_3 d$$
for some $\varepsilon_1, \varepsilon_2, \varepsilon_3\in\{-1,1\}$.
Now, it is easily checked that
\begin{equation}
\Phi(\left(\begin{array}{cc} a & c+id \\ c-id &b
\end{array}\right))\in \left(\begin{array}{cc} \varepsilon_1a & \varepsilon_2c+i\varepsilon_3d \\\varepsilon_2c+i\varepsilon_3d &
\varepsilon_1b
\end{array}\right)+{\mathbb R}I_2
\end{equation}
for some $\varepsilon_1,\varepsilon_2,\varepsilon_3\in\{-1,1\}$, and
the Claim 2 is true.

Replacing $\Phi$ by $\varepsilon_3(\Phi-f)$ if necessary, by Claim
2, we may assume that $\varepsilon_3\equiv 1$ and
\begin{equation} \Phi(A)=\Phi(\left(\begin{array}{cc} a & c+id \\
c-id &b
\end{array}\right))=\left(\begin{array}{cc} \varepsilon_1(A)a & \varepsilon_2(A)c+i d \\
\varepsilon_2(A)c-id & \varepsilon_1(A)b
\end{array}\right)\end{equation}
 for every $A\in {\bf H}_2$.

To determine the sign functions $\varepsilon_1,\varepsilon_2$ it is
enough to consider their behaviors on ${\bf H}_2^0$.

Let ${\mathcal M}=\{A\in{\bf H}_2^0:
\varepsilon_1(A)=\varepsilon_2(A)\}$ and ${\mathcal N}=\{B\in{\bf
H}_2^0 : \varepsilon_1(B)\not=\varepsilon_2(B)\}$.

{\bf Claim 2.} Either ${\mathcal M}={\bf H}_2^0$ or ${\mathcal
N}={\bf H}_2^0$.

   For
any $A=\left(\begin{array}{cc} a & c+id \\
c-id &-a
\end{array}\right),B=\left(\begin{array}{cc} b &e+if \\ e-if &
-b\end{array}\right)\in{\bf H}_2^0$, writing
$\varepsilon_1=\varepsilon_j(A)$ and $\eta_j=\varepsilon_j(B)$, a
simple computation shows that
$$AB-BA=2\left(\begin{array}{cc} i(de-cf) & ae-bc+i(af-bd) \\
-ae+bc+i(af-bd) &-i(de-cf)
\end{array}\right)$$ and
$$\begin{array}{rl} & \Phi(A)\Phi(B)-\Phi(B)\Phi(A)\\ =&2\left(\begin{array}{cc} i(\eta_2de-\varepsilon_2cf) & \varepsilon_1\eta_2ae-\varepsilon_2\eta_1bc+i(\varepsilon_1af-\eta_1bd) \\
-\varepsilon_1\eta_2ae+\varepsilon_2\eta_1bc+i(\varepsilon_1af-\eta_1bd)
&-i(\eta_2de-\varepsilon_2cf)
\end{array}\right). \end{array}$$
Since $w(\Phi(A)\Phi(B)-\Phi(B)\Phi(A))=w(AB-BA)$, one gets
$$\begin{array}{rl} &(\eta_2de-\varepsilon_2cf)^2+(\varepsilon_1\eta_2ae-\varepsilon_2\eta_1bc)^2+(\varepsilon_1af-\eta_1bd)^2\\
=& (de-cf)^2+ (ae-bc)^2+ (af-bd)^2,\end{array}
$$
that is,
$$\begin{array}{rl} & d^2e^2+c^2f^2+a^2f^2+b^2d^2-2df(\varepsilon_2\eta_2ce+\varepsilon_1\eta_1 ab)-2\varepsilon_1\varepsilon_2\eta_1\eta_2
abce \\ = & d^2e^2+c^2f^2+a^2f^2+b^2d^2-2df(ce+ab)-2abce,
\end{array}
$$
which gives \begin{equation}
df(\varepsilon_2\eta_2ce+\varepsilon_1\eta_1 ab)+
\varepsilon_1\varepsilon_2\eta_1\eta_2 abce  =  df(ce+ab)+abce.
\end{equation}

Assume  that both ${\mathcal M}$ and ${\mathcal N}$ are not empty.
Obviously, we can require $\varepsilon_1(A)=\varepsilon_2(A)$ if $ac=0$. So, ${\mathcal Q}=\{ A=\left(\begin{array}{cc} a & c+id \\
c-id &-a
\end{array}\right)\in{\bf H}_2^0:
ac=0\}\subseteq {\mathcal M}\cap{\mathcal N}$.  If one of ${\mathcal
M}$ and ${\mathcal N}$ is a subset of ${\mathcal Q}$, then the claim
is true. Assume that none of ${\mathcal M}$ and ${\mathcal N}$ is a
subset of ${\mathcal Q}$. We show that this leads to a
contradiction.

Let ${\mathcal N}_1={\mathcal N}\setminus{\mathcal Q}$. Then
${\mathcal
N}_1$ is not a empty set, and $B=\left(\begin{array}{cc} b &e+if \\
e-if & -b\end{array}\right)\in{\mathcal N}_1$ implies that
$be\not=0$. For
any $A=\left(\begin{array}{cc} a & c+id \\
c-id &-a
\end{array}\right)\in{\mathcal M}$, $B=\left(\begin{array}{cc} b &e+if \\ e-if &
-b\end{array}\right)\in{\mathcal N}$, since
$\varepsilon_2=\varepsilon_1=\varepsilon\in\{-1,1\}$ and
$\eta_2=-\eta_1=\eta\in\{-1,1\}$, Eq.(4.9) gives
$$ df\varepsilon\eta( ab-ce)-abce=df(ab+ce)+abce.
$$
Thus,

  if $\varepsilon\eta=1$, one gets $dfce=-abce$, that is,
$dfc=-abc$ as $be\not=0$;

 if $\varepsilon\eta=-1$, one gets
$dfab=-abce$, that is, $adf=-ace$ as $be\not=0$.

Assume that $f=0$ for some $\left(\begin{array}{cc} b &e+if \\ e-if
& -b\end{array}\right)\in{\mathcal N}_1 $; then we must have $ac=0$
for
all $\left(\begin{array}{cc} a & c+id \\
c-id &-a
\end{array}\right)\in{\mathcal M}$, which is a
contradiction. Thus, for all $B=\left(\begin{array}{cc} b &e+if \\
e-if & -b\end{array}\right)\in{\mathcal N}_1$, we have $bef\not=0$.
Hence,    for any $ A\in{\mathcal M}, B\in{\mathcal N}_1$,
$$\varepsilon(A)\varepsilon_1(B)=1\ {\rm and}\ c\not=0\Rightarrow df=-ab;$$
$$\varepsilon(A)\varepsilon_1(B)=-1 \ {\rm and}\ a\not=0\Rightarrow df=-ce.$$ Fix
some $A,B$ as above. Take $D=\left(\begin{array}{cc} x &y+iz \\ y-iz
& -x\end{array}\right)\in{\bf H}_2^0$ so that $xyz\not=0$,
$\frac{z}{x}\not\in\{\frac{d}{a},\frac{f}{b}\}$ and
$\frac{z}{y}\not\in\{ \frac{d}{c}, \frac{f}{e}\}$. Then it is easily
checked that $D\not ={\mathcal M}\cup{\mathcal N}={\bf H}_2^0$, a
contradiction. So, we must have ${\mathcal M}={\bf H}_2^0$ or
${\mathcal N}={\bf H}_2^0$.

{\bf Claim 3.} If ${\mathcal M}={\bf H}_2^0$, then $\Phi $ has the
form ($1^\circ$) or (2$^\circ$).

Let ${\mathcal M}_+=\{B\in{\mathcal M}: \varepsilon_1(B)=1\}$ and
${\mathcal M}_-=\{B\in{\mathcal M}: \varepsilon_1(B)=-1\}$. Then
 ${\bf H}_2^0={\mathcal M}={\mathcal M}_+\cup
{\mathcal M}_-$ and ${\mathcal M}_+\cap {\mathcal
M}_-=\{\left(\begin{array}{cc} 0 & if \\  -if & 0\end{array}\right):
f\in{\mathbb R}\}$. It is clear that $\Phi(A)=A$ if $A\in {\mathcal
M}_+$ and $\Phi(A)=-A^t$ if $A\in{\mathcal M}_-$.

For any $A=\left(\begin{array}{cc} a & c+id \\  c-id &
-a\end{array}\right)\in{\mathcal M}_+$ and $B=\left(\begin{array}{cc} b & e+ if \\
e-if & -b\end{array}\right)\in{\mathcal M}_-$, by Eq.(4.9) we have
$$df(ab+ce)=0.$$

Assume $df=0$; then the above equation is always true. If $f=0$,
then $B$ is a real matrix and $\Phi(B)=-B^t=-B$. Letting $h(B)$
absorb a $-1$ we may require that $B\in{\mathcal M}_+$. Similarly,
if $d=0$, we may rearrange if necessary so that $A\in{\mathcal
M}_-$. Hence we may require that one of ${\mathcal M}_\pm$ contains
no real matrices.

If one of ${\mathcal M}_\pm$ consists of real matrices, we already
prove that $\Phi$ has the form ($1^\circ$) or ($2^\circ$)

 Assume that   ${\mathcal M}_+$ and ${\mathcal M}_-$ contain respectively   non-real matrices $A$ and $B$; then $df\not=0$. It follows that
$$ab+ce=0.$$
If $abce\not=0$, we get
$$\frac{e}{b}=-\frac{a}{c}.$$
Take $D=\left(\begin{array}{cc} x &y+iz \\ y-iz &
-x\end{array}\right)\in{\bf H}_2^0$ with $xyz\not=0$,
$\frac{y}{x}\not\in\{\frac{c}{a},\frac{e}{b}\}$. Then either
$D\in{\mathcal M}_+$ or $D\in{\mathcal M}_-$. However,
$D\in{\mathcal M}_+$ implies that
$\frac{y}{x}=-\frac{b}{e}=\frac{c}{a}$ and $D\in{\mathcal M}_-$
implies that $\frac{y}{x}=-\frac{a}{c}=\frac{e}{b}$, contradicting
to the choice of $D$. Hence we always have $abce=0$, that is, at
lest one of $a,b,c,e$ is zero. Without loss of generality, assume
that $ac\not=0$; then $be=0$. In fact  we have $b=e=0$ since
$ab+ce=0$. This forces that ${\mathcal
M}_-=\{\left(\begin{array}{cc} 0& if \\   -if & 0\end{array}\right)
: f\in{\mathbb R}\}.$  and therefore ${\mathcal M}_+={\bf H}_2^0$.
In this case we have $\Phi(A)=A$ for all $A\in{\bf H}_2^0$ and
$\Phi$ has the form ($1^\circ$). If $be\not=0$ and $ac=0$, one gets
$a=c=0$ and thus
$${\mathcal M}_+\subseteq {\mathcal
R}=\{\left(\begin{array}{cc} u & w+iv \\  w-iv &
-u\end{array}\right) : v=0 \ \mbox{\rm or } u=w=0\}.$$ So we may
require that ${\mathcal M}_-={\bf H}_2^0$ and $\Phi(A)=-A^t$ for
every $A\in{\bf H}_2^0$, which implies that $\Phi$ has the form
($2^\circ$). If   $ac=be=0$ for any $A,B$ with $df\not=0$, then we
get a contradiction that $D=\left(\begin{array}{cc} x &y+iz \\ y-iz
& -x\end{array}\right)\in{\bf H}_2^0$ with $xyz\not=0$ does not in
${\mathcal M}_+\cup{\mathcal M}_-={\bf H}_2^0$. This completes the
proof of Claim 3.

{\bf Claim 4.} If ${\mathcal N}={\bf H}_2^0$, then $\Phi $ has the
form ($3^\circ$) or ($4^\circ$).

Let ${\mathcal N}_+=\{B\in{\mathcal N}: \varepsilon_1(B)=1\}$ and
${\mathcal N}_-=\{B\in{\mathcal N}: \varepsilon_1(B)=-1\}$. Then
 ${\bf H}_2^0={\mathcal N}={\mathcal N}_+\cup
{\mathcal N}_-$ and still, ${\mathcal N}_+\cap {\mathcal
N}_-=\{\left(\begin{array}{cc} 0 & if \\
-if & 0\end{array}\right): f\in{\mathbb R} \}$. Clearly,
$\Phi(A)=\Psi(A)$ if $A\in{\mathcal N}_+$ and $\Phi(A)=-\Psi(A)^t$
if $A\in{\mathcal N}_-$.

Note that,   for any $B_1,B_2\in{\mathcal N}_+$ or
$B_1,B_2\in{\mathcal N}_-$ we have
$w([B_1,B_2])=w([\Phi(B_1),\Phi(B_2)])$ by Eq.(4.9). Also, if $B$ is
real, then $\Phi(B)=-B$. Thus, with no loss of generality we may
assume that all real matrices are contained in ${\mathcal N}_+$.

For any $A=\left(\begin{array}{cc} a & c+id \\  c-id &
-a\end{array}\right)\in{\mathcal N}_+$ and $B=\left(\begin{array}{cc} b & e+ if \\
e-if & -b\end{array}\right)\in{\mathcal N}_-$, by Eq.(4.9) we still
have
$$df(ab+ce)=0.$$
If for any $A\in{\mathcal N}_+$ and $B\in{\mathcal N}_-$ we always
have $df=0$ whenever $(a,c)\not=(0,0)$, then we must have $d=0$ for
any $A\in{\mathcal N}_+$, which means that ${\mathcal N}_+\subseteq
{\mathcal R} $. It is easily checked in this case that $\Phi$ has
the form ($3^\circ$). So, we may assume that $df\not=0$ for some $A$
with $(a,c)\not=(0,0)$ and $B$. It follows that $ab+ce=0$. The same
reason as that in Claim 3 reveals that $abce\not=0$ will lead to a
contradiction. Thus we must have $abce=0$. Since there exists
$A\in{\mathcal N}_+$ with $acd\not=0$ or $B\in{\mathcal N}_-$ with
$bef\not=0$, a similar argument as that in Claim 3 shows that the
prior case implies that ${\mathcal N}_-=\{\left(\begin{array}{cc} 0 & if \\
-if & 0\end{array}\right): f\in{\mathbb R} \}$ and hence $\Phi$ has
the form ($3^\circ$); the later case implies that ${\mathcal
N}_+={\mathcal R}$ and hence $\Phi$ has the form ($4^\circ$).

 \hfill$\Box$

\end{document}